\documentclass[11pt]{amsart}
\usepackage{pdfsync}
\usepackage{fullpage}
\usepackage{amsmath}
\usepackage{amssymb}
\usepackage{MnSymbol}
\usepackage{verbatim}
\usepackage[linktocpage]{hyperref}

\newtheorem{df}{Definition}[section]
\newtheorem{thm}[df]{Theorem}

\newtheorem{prop}[df]{Proposition}
\newtheorem{rem}[df]{Remark}

\newtheorem{lem}[df]{Lemma}

\newtheorem{cor}[df]{Corollary}
\newtheorem{conj}[df]{Conjecture}

\newcommand{\pf}{\textit{Proof.} }
\newcommand{\eq}{eqnarray*}

\newcommand{\lig}{\mathfrak{g}}
\newcommand{\lik}{\mathfrak{k}}
\newcommand{\ga}{\gamma}
\newcommand{\ti}{\tilde}

\title{On cyclic Higgs bundles}

\author{Song Dai\textsuperscript{1}}

\address{Song Dai\\
Center for Applied Mathematics of Tianjin University\\
Tianjin University\\
No.92 Weijinlu Nankai District\\
Tianjin\\
P.R.China 300072\\}

\email{song.dai@tju.edu.cn}

\author{Qiongling Li\textsuperscript{2}}

\address{Qiongling Li\\
Centre for Quantum Geometry of Moduli Spaces (QGM)\\
Department of Mathematics, Aarhus University\\
Ny Munkegade 118, Bldg. 1530\\
8000 Aarhus C\\
Denmark\\}

\address{Department of Mathematics\\
California Institute of Technology\\
1200 East California Boulevard\\
Pasadena, CA 91125\\}

\email{qiongling.li@gmail.com}

\date{}

\begin{document}
\maketitle

\begin{abstract}
In this paper, we derive a maximum principle for a type of elliptic systems and apply it to analyze the Hitchin equation for cyclic Higgs bundles. We show several domination results on the pullback metric of the (possibly branched) minimal immersion $f$ associated to cyclic Higgs bundles. Also, we obtain a lower and upper bound of the extrinsic curvature of the image of $f$. As an application, we give a complete picture for maximal $Sp(4,\mathbb{R})$-representations in the $2g-3$ Gothen components and the Hitchin components.
\end{abstract}

\footnotetext[1]{The author is supported by NSFC grant No. 11601369.}

\footnotetext[2]{Corresponding author, supported in part by the center of excellence grant `Center for Quantum Geometry of Moduli Spaces' from the Danish National Research Foundation (DNRF95).}
\section{Introduction}
Let $S$ be a closed, oriented surface of genus $g\geq 2$ and $G$ be a reductive Lie group. Let $\Sigma$ be a Riemann surface over $S$ and denote its canonical line bundle by $K_{\Sigma}$. A $G$-Higgs bundle over $\Sigma$ is a pair $(E,\phi)$ where $E$ is a holomorphic vector bundle and $\phi$ is a holomorphic section of $End(E)\otimes K_{\Sigma}$ plus extra condition depending on $G$. The non-abelian Hodge theory developed by Corlette \cite{Corlette}, Donaldson \cite{Donaldson}, Hitchin \cite{Hitchin87} and Simpson \cite{Simpson88}, provides a one-to-one correspondence between the moduli space of representations from $\pi_1(S)$ to $G$ with the moduli space of $G$-Higgs bundles over $\Sigma$. The correspondence is through looking for an equivariant harmonic map from $\widetilde{\Sigma}$ to the symmetric space $G/K$, where $K$ is the maximal compact subgroup of $G$, for a given representation $\rho$ or a given Higgs bundle $(E,\phi)$.

In this paper, we are interested in the direction of the non-abelian Hodge correspondence from the moduli space of Higgs bundles to the space of equivariant harmonic maps. More explicitly, given a polystable $G$-Higgs bundle $(E,\phi)$ on $\Sigma$, there exists a unique Hermitian metric $h$ compatible with $G$-structure satisfying the Hitchin equation
\begin{\eq}
F^{\nabla^{h}}+[\phi,\phi^{*_{h}}]=0,
\end{\eq}called the harmonic metric, which gives the equivariant harmonic map from $\widetilde{\Sigma}$ to $G/K$. So for a given Higgs bundle $(E,\phi)$, we would like to deduce geometric properties of the corresponding equivariant harmonic map: $\widetilde{\Sigma}\rightarrow G/K$.

We are particularly interested in the following $SL(n,\mathbb{C})$-Higgs bundles
\begin{\eq}
E=L_1\oplus L_2\oplus \cdots\oplus L_n, \quad \phi=
\left(
\begin{array}{cccccc}
0 & & & \gamma_n\\
\gamma_1 & 0 & & \\
&\ddots & \ddots & \\
& & \gamma_{n-1} & 0
\end{array}
\right):E\rightarrow E\otimes K_{\Sigma},
\end{\eq}
where 
$L_k$ is a holomorphic line bundle and $\gamma_k$ is a holomorphic section of $L_{k}^{-1}L_{k+1}K_{\Sigma}$, $k=1,\cdots, n$ $(L_{n+1}=L_1$). Suppose $\det{E}=\mathcal{O}$ and $\ga_k\neq 0$, $k=1,\cdots ,n-1$. 
Call such a Higgs bundle $(E,\phi)$ a cyclic Higgs bundle parameterized by $(\ga_1,\ga_2,\cdots,\ga_n)$. For $G$ a subgroup of $SL(n,\mathbb{C})$, we call $(E,\phi)$ a cyclic $G$-Higgs bundle if it is a $G$-Higgs bundle and it is cyclic as a $SL(n,\mathbb{C})$-Higgs bundle. 

The terminology ``cyclic Higgs bundles" first appeared in \cite{Bar}. Note that the notion here is a bit different from the one in \cite{Bar}, where the notion ``cyclic" there is referred to the group $G$. 
One may also view cyclic Higgs bundles as a special type of quiver bundles in \cite{quiver}. Cyclic Higgs bundles are special in $G$-Higgs bundles for $G$ of higher rank. The harmonic metric for a cyclic Higgs bundle is diagonal, making it possible to analyze the solution to the Hitchin equation and hence the corresponding harmonic map. So studying cyclic Higgs bundles could give us hint on predicting what may happen to general Higgs bundles.

If a representation $\rho:\pi_1(S)\rightarrow SL(n,\mathbb{C})$ does not correspond to a cyclic Higgs bundle over one Riemann surface $\Sigma$, it is still possible that $\rho$ corresponds to a cyclic Higgs bundle over another Riemann surface $\Sigma'$. 
By Labourie \cite{LabourieCyclic}, any Hitchin representation for $SL(n,\mathbb{R})$ cannot correspond to a cyclic Higgs bundle over a deformation family of Riemann surfaces and later Collier \cite{Collier} generalizes to a more general family of cyclic Higgs bundles.


If $n\geq 3$, 
the associated harmonic map for a cyclic Higgs bundle is conformal and hence is a (possibly branched) minimal immersion. 
In \cite{DaiLi}, the authors studied the pullback metric and curvature of the minimal immersion for cyclic Higgs bundles in the Hitchin component (the $q_n$ case). In this paper, we derive a maximum principle for the elliptic systems. The maximum principle is very useful for the Toda-type equation with function coefficient, which appears in the Hitchin equation for cyclic Higgs bundles. With this powerful tool, we generalize and improve the results in \cite{DaiLi} and discover some new phenomena.

\subsection{Monotonicity of pullback metrics}
Let $(E,\phi)$ be a cyclic $SL(n,\mathbb{C})$-Higgs bundle parameterized by $(\ga_1,\cdots,\ga_n)$, $n\geq 3$. Let $f$ be the corresponding harmonic map and it is in fact branched minimal. The Riemannian metric on $SL(n,\mathbb{C})/SU(n)$ is induced by the Killing form on $sl(n,\mathbb{C})$.
Then the pullback metric of $f$ is given by
\begin{\eq}
g=2n\text{tr}(\phi\phi^{*_{h}})dz\otimes d\bar{z},
\end{\eq} where $h$ is the harmonic metric.
Though at branch points $g=0$, we still call $g$ a ``metric".

There is a nature $\mathbb{C}^*$-action on the moduli space $M_{Higgs}$ of $SL(n,\mathbb{C})$-Higgs bundles:
\begin{\eq}
\mathbb{C}^*\times \mathcal{M}_{Higgs}&\longrightarrow & \mathcal{M}_{Higgs}\\
t\cdot (E,\phi)&=& (E,t\phi)\end{\eq}

\begin{thm}Let $(E,\phi)$ be a cyclic Higgs bundle. Then along the $\mathbb{C}^*$-orbit of $(E,\phi)$, outside the branched points, as $|t|$ increases, the pullback metric $g^t$ of the corresponding branched minimal immersions strictly increases.\end{thm}

If we integrate the pullback metric, it gives the Morse function $f$ (up to a constant scalar) on the moduli space of Higgs bundles as the $L^2$-norm of $\phi$: \begin{\eq} f(E,\phi)=\int_{\Sigma}\text{tr}(\phi\phi^*) \sqrt{-1}dz\wedge d\bar z.\end{\eq}
\begin{cor}\label{Monotonicity}
Let $(E,\phi)$ be a cyclic Higgs bundle. Then along the $\mathbb{C}^*$-orbit of $(E,\phi)$, the Morse function $f(E,t\phi)$ strictly increases as $|t|$ increases.
\end{cor}
\begin{rem}
The Morse function is the main tool to determine the topology of the moduli space of Higgs bundles, for example, in Hitchin \cite{Hitchin87, Hitchin92}, Gothen \cite{Gothen}. The monotonicity in Corollary \ref{Monotonicity} is not new. In fact, Hitchin in \cite{Hitchin87} showed that with respect to the K\"ahler metric on the moduli space, the gradient flow of the Morse function is exactly the $\mathbb{R}^*$-part of $\mathbb{C}^*$-action. Hence, along $\mathbb{C}^*$-orbit of any Higgs bundles $(E,\phi)$, the Morse function $f(E,t\phi)$ strictly increases as $|t|$ increases.
Here we improve the integral monotonicity to pointwise monotonicity along $\mathbb{C}^*$-orbit of cyclic Higgs bundles.
\end{rem}

Consider the family of cyclic Higgs bundles $(E,\phi^t)$ parameterized by $(\ga_1,\cdots,t\ga_n)$. For $t\in \mathbb{C}^*$, the family $(E,\phi^t)$ is gauge equivalent to $t^{\frac{1}{n}}\cdot(E,\phi)=(E,t^{\frac{1}{n}}\phi)$. 
If the cyclic Higgs bundle parameterized by $(\ga_1,\cdots,\ga_{n-1},0)$ is again stable, in this case $\sum_{i=1}^k\text{deg}(L_{n+1-k})<0$ for all $1\leq k\leq n-1$, we extend the monotonicity of $\mathbb{C}^*$-family to the $\mathbb{C}$-family.
\begin{thm}\label{domination}
Let $(E,\phi^t)$ be a cyclic Higgs bundle parameterized by $(\ga_1,\cdots,t\ga_n)$ for $t\in \mathbb{C}$. If the cyclic Higgs bundle parameterized by $(\ga_1,\cdots,\ga_{n-1},0)$ is stable, then outside the branched points, as $|t|$ increases, the pullback metric $g^t$ of corresponding branched minimal immersions for $(E,\phi^t)$ strictly increases. So does the Morse function.
\end{thm}
\begin{rem}
If the cyclic Higgs bundle parameterized by $(\ga_1,\cdots,\ga_{n-1},0)$ is stable, it lies in the moduli space of Higgs bundles and is fixed by the $\mathbb{C}^*$-action. Note that the fixed points of $\mathbb{C}^*$-action are exactly the critical points of the Morse function as shown in Hitchin \cite{Hitchin87}.
\end{rem}

\subsection{Curvature of cyclic Higgs bundles in the Hitchin component} By Hitchin's description \cite{Hitchin92} of the Higgs bundle in the Hitchin component, the cyclic Higgs bundles in the Hitchin component are of the following form
\begin{\eq}
E=K^{\frac{n-1}{2}}\oplus K^{\frac{n-3}{2}}\oplus\cdots\oplus K^{\frac{3-n}{2}}\oplus K^{\frac{1-n}{2}},\quad
\phi=
\left(
\begin{array}{cccccc}
0 & & & q_n\\
1 & 0 & &\\
& \ddots & \ddots &\\
& & 1 & 0
\end{array}
\right),
\end{\eq}
where $q_n$ is a holomorphic $n$-differential. We call such a Higgs bundle $(E,\phi)$ a cyclic Higgs bundle in the Hitchin component parameterized by $q_n$. If $q_n=0$, the Higgs bundle is called Fuchsian.

The corresponding harmonic map $f:\widetilde{\Sigma}\rightarrow SL(n,\mathbb{R})/SO(n)$ is a minimal immersion for $n\geq 3$. We want to investigate that, as an immersed submanifold, how the image $f(\widetilde{\Sigma})$ sits inside the symmetric space $SL(n,\mathbb{R})/SO(n)$.

\begin{thm}\label{Introcurvature} Let $(E,\phi)$ be a cyclic Higgs bundle in the Hitchin component parameterized by $q_n$. Let $\sigma$ be the tangent plane of the image of $f$, then the curvature $K_{\sigma}$ in $SL(n,\mathbb{R})/SO(n)$ satisfies $$-\frac{1}{n(n-1)^2}\leq K_{\sigma}<0.$$
\end{thm}
\begin{rem} The sectional curvature $K$ of $SL(n,\mathbb{R})/SO(n)$ and $SL(n,\mathbb{C})/SU(n)$ satisfies $-\frac{1}{n}\leq K\leq 0$ (see Proposition \ref{CurvatureSym}).
For general Higgs bundles, one should not expect there is such a nontrivial lower bound at immersed points. For example, in the case of cyclic Higgs bundles parametrized by $(\ga_1,\ga_2,\cdots,\ga_n)$, if $n-1$ terms of $\ga_i$'s have a common zero point, then the curvature of the tangent plane $\sigma$ at that point achieves the most negative, i.e., $K_{\sigma}=-\frac{1}{n}$.
\end{rem}
\begin{rem}
(1) The upper bound is shown in \cite{DaiLi}. Here we give a new proof. As shown in \cite{CL14}, along the family of Higgs bundles parameterized by $tq_n$, $K_{\sigma}^t$ approaches to $0$ away from the zeros of $q_n$ as $|t|\rightarrow\infty$. So the upper bound $K_{\sigma}<0$ is sharp.\\
(2) The lower bound $-\frac{1}{n(n-1)^2}$ can only be achieved at some point in the case $n=2,3$.\\
(3) In the Fuchsian case, i.e. $q_n=0$, the sectional curvature $K_{\sigma}$ is $-\frac{6}{n^2(n^2-1)}$. However, it is strictly larger than the lower bound of $K_{\sigma}$ for $q_n\neq 0$ case when $n>3$. Hence, one cannot expect the curvature in Fuchsian case could serve as a lower bound of $K_{\sigma}$ for general Hitchin representations.
\end{rem}
For details one may see the remarks in the end of Section \ref{curvature}.

\subsection{Comparison inside the real Hitchin fibers}
Fix a Riemann surface $\Sigma$, the Hitchin fibration is a map from the moduli space of $SL(n,\mathbb{C})$-Higgs bundles over $\Sigma$ to the direct sum of the holomorphic differentials
\begin{\eq}
h:M_{Higgs} &\longrightarrow&\bigoplus\limits_{j=2}^nH^0(\Sigma, K^j)\ni (q_2,q_3,\cdots,q_n).\\
(E,\phi)&\mapsto& (\text{tr}(\phi^2),\text{tr}(\phi^3),\cdots,\text{tr}(\phi^n))\end{\eq}
Note that cyclic Higgs bundles $(E,\phi)$ lie in the Hitchin fiber at $(0,\cdots,0,n\cdot q_n)$, where $q_n=(-1)^{n-1}\det(\phi)$. There is one special point in each Hitchin fiber at $(0,\cdots,0,n\cdot q_n)$: the cyclic Higgs bundle in the Hitchin component parametrized by $q_n$.

In Proposition \ref{harmonicmetriccomparison}, we show that the harmonic metric in the cyclic Higgs bundle in the Hitchin component dominates the ones for other cyclic $SL(n,\mathbb{R})$-Higgs bundles in the same Hitchin fiber in a certain sense.

As the applications in lower rank $n=2,3,4$, we compare the pullback metric of the harmonic map for the cyclic Higgs bundle in the Hitchin component with the ones for other cyclic $SL(n,\mathbb{R})$-Higgs bundles in the same Hitchin fiber at $(0,\cdots, 0,n\cdot q_n)$.
\begin{thm}\label{pullback}
Let $(\tilde{E},\tilde{\phi})$ be a cyclic Higgs bundle in the Hitchin component parameterized by $q_n$ and $(E,\phi)$ be a distinct cyclic $SL(n,\mathbb{R})$-Higgs bundle in Section \ref{s2} such that $\det\phi=(-1)^{n-1}q_n$. For $n=2,3,4$, the pullback metrics $g,\ti{g}$ of the corresponding harmonic maps satisfy $g<\ti{g}.$
\end{thm}

Under the assumptions above, the Morse function satisfies $f(E,\phi)<f(\ti E,\ti\phi)$.

By Hitchin's work in \cite{Hitchin87}, all polystable $SL(2,\mathbb{R})$-Higgs bundles with nonvanishing Higgs field are cyclic. We can then directly apply Theorem \ref{pullback} to $SL(2,\mathbb{R})$-representations, we recover the following result shown in \cite{DominationFuchsian}.
\begin{cor}For any non-Fuchsian reductive $SL(2,\mathbb{R})$-representation $\rho$ and any Riemann surface $\Sigma$, there exists a Fuchsian representation $j$ such that the pullback metric of the corresponding $j$-equivariant harmonic map $f_j:\widetilde{\Sigma}\rightarrow \mathbb{H}^2$ dominates the one for $f_{\rho}$. \end{cor}
\begin{rem}
Deroin and Tholozan in \cite{DominationFuchsian} show a stronger result by comparing Fuchsian representations with all $SL(2,\mathbb{C})$-representations and the condition being reductive can be removed by separate consideration. Inspired by this result, they conjecture that in the Hitchin fiber, the Hitchin section maximizes the translation length. Our Theorem \ref{pullback} here is exactly in the same spirit, but using the pullback metric rather than the translation length.
\end{rem}

We expect that Theorem \ref{pullback} holds for general Higgs bundles rather than just cyclic Higgs bundles.
\begin{conj}\label{main}
Let $(\tilde{E},\tilde{\phi})$ be a Higgs bundle in the Hitchin component and $(E,\phi)$ be a distinct $SL(n,\mathbb{R})$-Higgs bundle in the same Hitchin fiber at $(q_2,q_3,\cdots,q_n)$. Then the pullback metrics $g,\ti{g}$ of corresponding harmonic maps satisfy $g<\ti{g}.$ As a result, the Morse function satisfies $f(E,\phi)<f(\ti E,\ti\phi)$.
\end{conj}
\subsection{Maximal $Sp(4,\mathbb{R})$-representations} For each reductive representation $\rho$ into a Hermitian Lie group $G$, we can define a Toledo integer $\tau(\rho)$ satisfying the Milnor-Wood inequality $|\tau(\rho)|\leq \text{rank}(G)(g-1)$. The representation $\rho$ with $|\tau(\rho)|=\text{rank}(G)(g-1)$ is called maximal. Maximal representations are Anosov \cite{Burger} and hence discrete and faithful.

In the case for $Sp(4,\mathbb{R})$, there are $3\cdot 2^{2g}+2g-4$ connected components of maximal representations containing $2^{2g}$ isomorphic components of Hitchin representations \cite{Hitchin92} and $2g-3$ exceptional components called Gothen components \cite{Gothen}.
Labourie in \cite{LabourieCyclic} shows that any Hitchin representation corresponds to a cyclic Higgs bundle in the Hitchin component parametrized by $q_4$ over a unique Riemann surface. Together with the description in \cite{Gothen,BradlowDeformation} and Collier's work \cite{Collier}, any maximal representation for $Sp(4,\mathbb{R})$ in the Gothen components corresponds to a cyclic Higgs bundle over a unique Riemann surface $\Sigma$ of the form
\begin{\eq}
E=N\oplus NK^{-1}\oplus N^{-1}K\oplus N^{-1}, \quad \phi=
\left(
\begin{array}{cccccc}
0 & && \nu\\
1 &0 & &\\
&\mu &0 & \\
& & 1 & 0
\end{array}
\right),
\end{\eq} where $g-1<\deg N<3g-3$. 
Note that if $N=K^{\frac{3}{2}}$, the above Higgs bundle corresponds to a Hitchin representation. As a result, for any $Sp(4,\mathbb{R})$-representation in the Hitchin components or Gothen components, there is a unique $\rho$-equivariant minimal immersion of $\widetilde{S}$ in $Sp(4,\mathbb{R})/U(2)$. Recently, this result is reproved and generalized to maximal $SO(2,n)$-representations in Collier-Tholozan-Toulisse \cite{CollierTholozan}.

For each Riemann surface, the above cyclic Higgs bundles with $\nu=0$ play a similar role as the Fuchsian case: they are the fixed points of the $\mathbb{C}^*$-action. We call the corresponding representations $\mu$-Fuchsian representations. The only difference with the Fuchsian case is that they form a subset inside each component rather than one single point since $\mu\in H^0(N^{-2}K^3)$ has many choices. As a corollary of Theorem \ref{domination}, the space of $\mu$-Fuchsian representations serves as the minimum set in its component of maximal $Sp(4,\mathbb{R})$ representations in the following sense.
\begin{cor}
For any maximal $Sp(4,\mathbb{R})$-representation $\rho$ in the $2g-3$ Gothen components (or the Hitchin components), there exists a $\mu$-Fuchsian (or Fuchsian) representation $j$ in the same component of $\rho$ such that the pullback metric of the unique $j$-equivariant minimal immersion $f_j:\widetilde{S}\rightarrow Sp(4,\mathbb{R})/U(2)$ is dominated by the one for $f_{\rho}$.
\end{cor}
To consider the curvature, as a corollary of Theorem \ref{Introcurvature}, we have
\begin{cor}\label{HitchinCurvatureIntro}
For any Hitchin representation $\rho$ for $Sp(4,\mathbb{R})$, the sectional curvature $K_{\sigma}$ in $Sp(4,\mathbb{R})/U(2)$ of the tangent plane $\sigma$ of the unique $\rho$-equivariant minimal immersion satisfies \\
(1) $K_{\sigma}=-\frac{1}{40}$, if $\rho$ is Fuchsian;\\
(2) $-\frac{1}{36}<K_{\sigma}<0$ and $\exists~ p$ such that $K_{\sigma}(p)<-\frac{1}{40},$ if $\rho$ is not Fuchsian.
\end{cor}

Similarly, we also obtain an upper and lower bound on the curvature of minimal immersions for maximal representations.
\begin{thm}\label{MaximalCurvatureIntro}
For any maximal representation $\rho$ for $Sp(4,\mathbb{R})$ in each Gothen ccomponent, the sectional curvature $K_{\sigma}$ in $Sp(4,\mathbb{R})/U(2)$ of tangent plane $\sigma$ of the uniuqe $\rho$-equivariant minimal immersion satisfies \\
(1) $-\frac{1}{8}\leq K_{\sigma}<-\frac{1}{40}$ and the lower bound is sharp, if $\rho$ is $\mu$-Fuchsian;\\
(2) $-\frac{1}{8}\leq K_{\sigma}<0$, if $\rho$ is not $\mu$-Fuchsian.
\end{thm}
\begin{rem} As shown in \cite{CL14},\cite{Mochizuki}, along the family of $(E,t\phi)$, away from zeros of $det(\phi)\neq 0$, the sectional curvature goes to zero as $|t|\rightarrow \infty$. So the upper bounds in Part (2) in both Corollary \ref{HitchinCurvatureIntro} and Theorem \ref{MaximalCurvatureIntro} are sharp.
The sectional curvature $K$ in $Sp(4,\mathbb{R})/U(2)$ satisfies $-\frac{1}{4}\leq K\leq 0$. So the lower bounds in Corollary \ref{HitchinCurvatureIntro} and Theorem \ref{MaximalCurvatureIntro} are nontrivial.
\end{rem}

As an immediate corollary of Theorem \ref{pullback} for $n=4$, comparing maximal representations in the Gothen components with Hitchin representations, we obtain
\begin{cor}
For any maximal $Sp(4,\mathbb{R})$-representation $\rho$ in the $2g-3$ Gothen components, there exists a Hitchin representation $j$ such that the pullback metric of the unique $j$-equivariant minimal immersion $f_j:\widetilde{S}\rightarrow Sp(4,\mathbb{R})/U(2)$ dominates the one for $f_{\rho}$.
\end{cor}

\subsection{Maximum principle}
We derive a maximum principle for the elliptic systems. It is the main tool we use throughout this paper.

Basically, we consider the following linear elliptic system
\begin{eqnarray*}
\triangle_{g} u_i+<X,\nabla u_i>+\sum_{j=1}^{n}c_{ij}u_j=f_i, \quad 1\leq i\leq n.
\end{eqnarray*}
Roughly speaking, suppose the functions $c_{ij}$ satisfy the following assumptions: \\
$(a)$ cooperative: $c_{ij}\geq 0,~ i\neq j$,\\
$(b)$ column diagonally dominant: $\sum_{i=1}^{n}c_{ij}\leq 0,~ 1\leq j\leq n$,\\
$(c)$ fully coupled: the index set $\{1,\cdots,n\}$ cannot be split up in two disjoint nonempty sets $\alpha,\beta$ such that $c_{ij}\equiv 0$ for $i\in\alpha,j\in \beta.$\\
Then the maximum principle holds, that is, if $f_i\leq 0$ for $1\leq i\leq n$, then $u_i\geq 0$ for $1\leq i\leq n$.
The precise statement is Lemma \ref{mp}.

In the literature, it is common to require there exists a positive supersolution, which is equivalent to the maximum principle, see \cite{LM}. So for function coefficients, people usually suppose $c_{ij}$ satisfy the row sum condition $\sum_{j=1}^{n}c_{ij}\leq 0,~ 1\leq i\leq n$, say \cite{PW}. The column sum condition $\sum_{i=1}^{n}c_{ij}\leq 0,~ 1\leq j\leq n$, or in other words column diagonally dominant condition, rarely appeared in the literature. The similar column sum condition first appeared in \cite{LM}, Theorem 3.3.

To the knowledge of the authors, the maximum principles in the literature seem not to directly imply our maximum principle Lemma \ref{mp}. We also remark that our proof is more elementary.

\subsection*{Structure of the article.} The article is organized as follows. In Section \ref{pre}, we recall some fundamental results about the Higgs bundle and introduce the cyclic Higgs bundles. In Section \ref{MaxP}, we show a maximum principle for the elliptic systems, the main tool of this article. In Section \ref{monotonicity}, we show the monotonicity of the pullback metrics of the branched minimal immersions. In Section \ref{curvature}, we find out a lower and upper bound for the extrinsic curvature of the minimal immersions for cyclic Higgs bundles in the Hitchin component. In Section \ref{comparison}, we compare the harmonic metrics of cyclic Higgs bundle in the Hitchin component with other cyclic $SL(n,\mathbb{R})$-Higgs bundles in the same Hitchin fiber. In Section \ref{maximal}, we apply our results to maximal $Sp(4,\mathbb{R})$-representations.

\subsection*{Acknowledgement} The authors wish to thank Nicolas Tholozan for suggesting the problem of looking for a lower bound for the extrinsic curvature of the harmonic map. The authors acknowledge support from U.S. National Science Foundation grants DMS 1107452, 1107263, 1107367 ``RNMS: GEometric structures And Representation varieties" (the GEAR Network).

\section{Preliminaries}\label{pre}
In this section, we recall some facts in the theory of the Higgs bundles. One may refer \cite{Bar}\cite{DaiLi}\cite{LabourieCyclic}. Let $\Sigma$ be a closed Riemann surface of genus $\geq 2$ and $K=K_{\Sigma}$ be the canonical line bundle over $\Sigma$. For $p\in\Sigma$, let $\pi_1=\pi_1(\Sigma,p)$ be the fundamental group of $\Sigma$. Let $\tilde{\Sigma}$ be the universal cover of $\Sigma$.

A $SL(n,\mathbb{C})$-Higgs bundle over $\Sigma$ is a pair $(E,\phi)$, where $E$ is a holomorphic vector bundle with $\det E=\mathcal{O}$ and $\phi$ is a trace-free holomorphic section of $End(E)\bigotimes K$. We call $(E,\phi)$ is stable if for any proper $\phi$-invariant holomorphic subbundle $F$,
$\frac{\text{deg}F}{\text{rank}F}<\frac{\text{deg}E}{\text{rank}E}$. We call $(E, \phi)$ is polystable if $(E,\phi)$ is a direct sum of stable Higgs bundles of degree $0$.
\subsection{Higgs bundles and harmonic maps}
\begin{thm}(Hitchin \cite{Hitchin87} and Simpson \cite{Simpson88})\label{dh}
Let $(E,\phi)$ be a stable $SL(n,\mathbb{C})$-Higgs bundle. Then there exists a unique Hermitian metric $h$ on $E$ compatible with $SL(n,\mathbb{C})$ structure, called the harmonic metric, solving the Hitchin equation
\begin{eqnarray}\label{hit eq}
F^{\nabla^h}+[\phi,\phi^{*_h}]=0, \label{curvature equation}
\end{eqnarray}
where ${\nabla^h}$ is the Chern connection of $h$, in local holomorphic trivialization,
\begin{\eq}
F^{\nabla^h}=\overline{\partial}(h^{-1}\partial h),
\end{\eq}
and $\phi^{*_{h}}$ is the adjoint of $\phi$ with respect to $h$, in the sense that
\begin{\eq}
h(\phi(u),v)=h(u,\phi^{*_{h}}(v))\in K,\quad u,v\in E
\end{\eq}
in local frame, $\phi^{*_{h}}=\bar{h}^{-1}\bar{\phi}^{\top}\bar{h}$.
\end{thm}
Denote\begin{\eq}
&&G=SL(n,\mathbb{C}),\quad K=SU(n)\\
&&\mathfrak{g}=sl(n,\mathbb{C}),\quad \mathfrak{k}=su(n),\quad \mathfrak{p}=\{X\in sl(n,\mathbb{C}):\bar X^t=X\},\quad\mathfrak{g}=\mathfrak{k}\oplus \mathfrak{p}.\\
&&\text{The Killing form } B(X,Y)=2n\cdot \text{tr}(XY)��
\end{\eq}
The harmonic metric $h$ gives rise to a flat $SL(n,\mathbb{C})$ connection $D=\nabla^h+\Phi=\nabla^{h}+\phi+\phi^{*_{h}}$. The holonomy of $D$ gives a representation $\rho:\pi_1\to SL(n,\mathbb{C})$ and the bundle $(E,D)$ is isomorphic to $\widetilde{\Sigma}\times_{\rho}\mathbb{C}^n$ with the associated flat connection. A Hermitian metric $h$ on $E$ is equivalent to a reduction $i:P_K\rightarrow P_G$ from unimodule frame bundle $P_G=\widetilde{\Sigma}\times_{\rho}G$ of $E=\widetilde{\Sigma}\times_{\rho}\mathbb{C}^n$ to the unitary frame bundle $P_K$ of $E$ with respect to $h$. Then it descends to be a section of $P_G/K=\widetilde{\Sigma}\times_{\rho}G/K$ over $\Sigma$. Equivalently, it gives a $\rho$-equivariant map $f: \tilde{\Sigma}\to G/K$. Denote the bundle $\tilde{P}_K$ be the pullback of the principle $K$-bundle $G\rightarrow G/K$ by $f$. Note that $\pi^*P_K=\tilde{P}_K$, where $\pi$ is the covering map $\pi: \tilde{\Sigma}\to \Sigma$. 
The Maurer-Cartan form $\omega$ of $G$ gives a flat connection on $P_G$, we still use $\omega$ to denote the connection. It coincides with the flat connection $D$. Consider $i^{*}\omega$, which is a $\lig$-value one form on $P_K$. Decomposing $i^{*}\omega=A+\Phi$ from $\lig=\lik\oplus\mathfrak{p}$, where $A$ is $\lik$-valued and $\Phi$ is $\mathfrak{p}$-valued. Then $A\in\Omega^1(P_K,\mathfrak{k})$ is a principal connection on $P_K$ and $\Phi$ is a section of $T^*\Sigma\otimes (P_{K}\times_{Ad_{K}}\mathfrak{p})$. By complexification, $\Phi$ is also a section of
\begin{\eq}
(T^*{\Sigma}\otimes\mathbb{C})\otimes ({P}_{K}\times_{Ad_{K}}\mathfrak{p}\otimes\mathbb{C})&=&(T^*{\Sigma}\otimes\mathbb{C})\otimes ({P}_{K^{\mathbb{C}}}\times_{Ad_{K^{\mathbb{C}}}}\mathfrak{p}^{\mathbb{C}})
\\
&=&(K\oplus\bar{K})\otimes ({P}_G\times_{Ad_G}\mathfrak{g})=(K\oplus\bar{K})\otimes End_0(E)\end{\eq}where $End_0(E)$ the trace-free endormorphism bundle of $E$. With respect to the decomposition $(K\oplus\bar{K})\otimes End_0(E)$, $\Phi=\phi+\phi^*$.

With respect to the decomposition $\mathfrak{g}=\mathfrak{k}+\mathfrak{p}$, we can decompose $\omega=\omega^{\mathfrak{k}}+\omega^{\mathfrak{p}}$, where $\omega^{\mathfrak{k}}\in \Omega^1(G,\mathfrak{k}), \omega^{\mathfrak{p}}\in \Omega^1(G,\mathfrak{p})$. Moreover, $\omega^{\mathfrak{p}}$ descends to be an element in $\Omega^1(G/K,G\times_{Ad_K}\mathfrak{p})$. In fact, using the Maurer-Cartan form $\omega^{\mathfrak{p}}\in \Omega^1(G/K,G\times_{Ad_K}\mathfrak{p})$ over $G/K$: $T(G/K)\cong G\times_{Ad_K}\mathfrak{p}$. 
Then $\mathfrak{g}=\mathfrak{k}\oplus\mathfrak{p}$ gives an $Ad_K$-invariant orthogonal decomposition and the Killing form $B$ on $\mathfrak{g}$ is positive on $\mathfrak{p}$. The Killing form $B$ induces a Riemannian metric $\tilde{B}$ on $G/K$: for two vectors $Y_1,Y_2\in T_p(G/K)$, $$\ti B(Y_1,Y_2)=B(\omega^{\mathfrak{p}}(Y_1),\omega^{\mathfrak{p}}(Y_2)).$$ Then $f^*\omega^{\mathfrak{p}}$ is a section of $T^*\tilde{\Sigma}\otimes (\tilde{P}_K\times_{Ad_K}\mathfrak{p})$ over $\tilde{\Sigma}$.

By comparing the two decomposition of the Maurer-Cartan form $\omega$, we obtain: $$f^*\omega^{\mathfrak{p}}=\pi^*\Phi.$$ So for every tangent vector $X\in T\widetilde\Sigma$, under the isomorphism by the Maurer-Cartan form $$\omega^{\mathfrak{p}}:T(G/K)\cong G\times_{Ad_K}\mathfrak{p},$$
we have
\begin{equation}\label{KeyFormula}
\omega^{\mathfrak{p}}(f_*(X))=f^*\omega^{\mathfrak{p}}(X)=\pi^*\Phi(X)=\Phi(\pi_*(X)).
\end{equation}

We consider the pullback metric $g$ on $\Sigma$, $g=\pi_{*}f^{*}\tilde{B}$. Since $f$ is $\rho$-equivariant and $\tilde{B}$ is $G$-invariant, $g$ is well defined. Then $\forall X,Y\in T\Sigma$, locally choose any lift $\tilde{X},\tilde{Y}\in T\tilde{\Sigma}$,
\begin{\eq}
g(X,Y)=\tilde{B}(f_*(\tilde{X}),f_*(\tilde{Y}))=B(\omega^{\mathfrak{p}}(f_*(\tilde{X})),\omega^{\mathfrak{p}}(f_*(\tilde{Y}))=B(\Phi(X),\Phi(Y)).
\end{\eq}

Later in the paper, we may ignore this covering map $\pi$ for short.
Then we have $$\text{Hopf}(f)=g^{2,0}=2n\text{tr}(\phi\phi),\quad g^{1,1}=2n\text{tr}(\phi\phi^{*_{h}})dz\otimes d\bar{z}.$$ If $\text{Hopf}(f)=0$, then as a section of $K\otimes\bar{K}$, the Hermitian metric is \begin{\eq}g=g^{1,1}=2n\text{tr}(\phi\phi^{*_{h}})dz\otimes d\bar z.\end{\eq} The associated Riemannian metric of $g$ is $g+\bar{g}$ on $\Sigma$, i.e., $2n\text{tr}(\phi\phi^{*_{h}})dz\cdot d\bar{z}$, where
$$dz\cdot d\bar{z}=dz\otimes d\bar{z}+d\bar{z}\otimes dz=2|dz|^2=2(dx^2+dy^2).$$

We focus on the cyclic Higgs bundles introduced below.
\subsection{Cyclic Higgs bundles}\label{s1}
A cyclic Higgs bundle is a $SL(n,\mathbb{C})$-Higgs bundle $(E,\phi)$ of the following form
\begin{\eq}
E=L_1\oplus L_2\oplus \cdots\oplus L_n, \quad \phi=
\left(
\begin{array}{cccccc}
0 & & & \gamma_n\\
\gamma_1 & 0 & & \\
& \ddots & \ddots &\\
& & \gamma_{n-1} & 0
\end{array}
\right),
\end{\eq}
where $L_k$ is a holomorphic line bundle over $\Sigma$ and $\gamma_k$ is a holomorphic section of $L_{k}^{-1}L_{k+1}K$, $k=1,\cdots, n$. The subscript is counted modulo $n$, i.e., $n+1\equiv 1$. Here $\det {E}=\mathcal{O}$ and $\ga_k\neq 0$, $k=1,\cdots ,n-1$. If $\ga_n\neq 0$, $(E,\phi)$ is automatically stable, which implies the existence of the solution to the Hitchin equation $(\ref{hit eq})$. If $\gamma_n=0$, $(E,\phi)$ stable in this case means $\sum_{i=1}^k\text{deg}(L_{n+1-i})<0$ for all $1\leq k\leq n-1$.

Following the proof in Baraglia \cite{Bar}, Collier \cite{Collier,BrianThesis}, the harmonic metric is diagonal for cyclic Higgs bundles. We include the proof here for completeness.
\begin{prop}
For a cyclic Higgs bundle $(E,\phi)$, the harmonic metric $h$ is diagonal, i.e. $$h=\text{diag}(h_1,h_2,\cdots, h_n)$$ where each $h_k$ is a Hermitian metric on $L_{k}$. %
\end{prop}
\begin{pf}
For $\omega=e^{\frac{2\pi i}{n}}$, consider the holomorphic $SL(n,\mathbb{C})$-gauge transformation $g_{\omega}$:
$$g_{\omega}=\begin{pmatrix}\omega^{\frac{n-1}{2}}&&&\\&\omega^{\frac{n-3}{2}}&&
\\&&\ddots&\\&&&\omega^{\frac{1-n}{2}}\end{pmatrix}:E\rightarrow E$$
It acts on the Higgs field $\phi$ as follows
$$g_{\omega}\cdot \phi=g_{\omega}\phi g_{\omega}^{-1}=\omega\cdot \phi$$

Then the metric $hg_{\omega}^{*_h}g_{\omega}$ is a solution to the Higgs bundle $(g_{\omega}^{-1}\bar\partial_{E}g_{\omega},g_{\omega}^{-1}\phi g_{\omega})=(\bar\partial_E,\omega^{-1}\cdot \phi)$. Since $U(1)$-action does not change the harmonic metric, $hg_{\omega}^{*_h}g_{\omega}$ is also the solution to the Higgs bundle $(\partial_E,\phi)$. Hence, by the uniqueness of harmonic metrics,
$$h=hg_{\omega}^{*_h}g_{\omega}.$$
Then $h$ splits as $(h_1,h_2,\cdots,h_n)$.
\end{pf}\qed

Denote $L\otimes\bar{L}=|L|^2$, then $h_k$ is a smooth section of $|L_k|^{-2}$. Chosen a local holomorphic frame, we abuse $\ga_k$ to denote the local coefficient function of the section $\ga_k$.
Then locally the Hitchin equation is
\begin{\eq}\triangle \log h_k+|\gamma_k|^2h_k^{-1}h_{k+1}-|\gamma_{k-1}|^2h_{k-1}^{-1}h_k=0,\quad k=1,\cdots, n,
\end{\eq}
where $\triangle=\partial_z\bar{\partial}_{z}$, $|\ga_k|^2=\ga_k\bar{\ga}_k$ as a local function.

If $n\geq 3$, the Hopf differential of the harmonic map $\text{Hopf}(f)=\text{tr}(\phi^2)=0$. And $f$ is immersed at $p$ if and only if $\phi(p)\neq 0$. At point $p$ where $\phi(p)=0$, $f$ is branched at $p$. 
Then outside the branch points, the harmonic map is conformal, then minimal. The pullback metric is given by $$g=2n\text{tr}(\phi\phi^{*_{h}})=2n(\sum_{k=1}^{n}|\gamma_k|^2h_k^{-1}h_{k+1})dz\otimes d\bar{z}.$$
\begin{rem}
For $n=2$, we consider the $(1,1)$ part of the pullback metric $g$ instead.
\end{rem}

\subsection{Cyclic $SL(n,\mathbb{R})$-Higgs bundles}\label{s2}
A $SL(n,\mathbb{R})$-Higgs bundle over $\Sigma$ is a triple $(E,\phi, Q)$, where $(E,\phi)$ is a $SL(n,\mathbb{C})$-Higgs bundle and $Q$ is a non-degenerate holomorphic quadratic form on $E$ such that $Q(\phi u,v)=Q(u,\phi v)$ for $u,v\in E$. Such $(E,\phi,Q)$ corresponds to a representation
\begin{\eq}
\rho:\pi_1\to SL(n,\mathbb{R})\hookrightarrow SL(n,\mathbb{C}).
\end{\eq}

Here we consider the holomorphic quadratic form
\begin{equation*}
Q=\left(
\begin{array}{cccccc}
&&& 1\\
&& 1&\\
& \udots &&\\
1& &&
\end{array}
\right):E\xrightarrow{\cong} E^*.\end{equation*}

For $n=2m$, the cyclic $SL(n,\mathbb{R})$-Higgs bundle is of the following form
\begin{\eq}
E=L_1\oplus\cdots\oplus L_m\oplus L_m^{-1}\oplus\cdots\oplus L_1^{-1},\quad \phi=
\left(
\begin{array}{cccccccc}
0 & & & & & & & \nu\\
\gamma_1 & \ddots & & & & & &\\
& \ddots & 0 & & & & &\\
& & \gamma_{m-1} & 0 & & & &\\
& & & \mu & 0 & & &\\
& & & & \ga_{m-1} & 0 & &\\
& & & & & \ddots & \ddots &\\
& & & & & & \gamma_{1} & 0
\end{array}
\right).
\end{\eq}
By the uniqueness of the solution, $h=\text{diag}(h_1,\cdots,h_m,h_{m}^{-1},\cdots,h_1^{-1})$. Locally, the Hitchin equation is
\begin{\eq}
\triangle \log h_1+|\gamma_1|^2h_1^{-1}h_{2}-|\nu|^2h_{1}^{2}&=&0,\\
\triangle \log h_k+|\gamma_k|^2h_k^{-1}h_{k+1}-|\gamma_{k-1}|^2h_{k-1}^{-1}h_k&=&0,\quad k=2,\cdots, m-1,\\
\triangle \log h_m+|\mu|^2h_m^{-2}-|\gamma_{m-1}|^2h_{m-1}^{-1}h_m&=&0.
\end{\eq}
The pullback metric is $g=2n(|\nu|^2h_{1}^{2}+|\mu|^2h_m^{-2}+2\sum_{k=1}^{m-1}|\gamma_k|^2h_k^{-1}h_{k+1})dz\otimes d\bar{z}.$
\\

For $n=2m+1$, the cyclic $SL(n,\mathbb{R})$-Higgs bundle is of the following form
\begin{\eq}
E=L_1\oplus\cdots\oplus L_m\oplus \mathcal{O}\oplus L_m^{-1}\oplus\cdots\oplus L_1^{-1},\quad \phi=
\left(
\begin{array}{ccccccccc}
0 & & & & & & & &\nu\\
\gamma_1 & \ddots & & & & & & &\\
& \ddots & 0 & & & & & &\\
& & \gamma_{m-1} & 0 & & & & &\\
& & & \mu & 0 & & & &\\
& & & & \mu & 0 & & &\\
& & & & & \ga_{m-1} & 0 & &\\
& & & & & & \ddots & \ddots &\\
& & & & & & & \gamma_{1} & 0
\end{array}
\right).
\end{\eq}
In this case, $h=\text{diag}(h_1,\cdots,h_m,1,h_{m}^{-1},\cdots,h_1^{-1})$. Locally, the Hitchin equation is
\begin{\eq}
\triangle \log h_1+|\gamma_1|^2h_1^{-1}h_{2}-|\nu|^2h_{1}^{2}&=&0,\\
\triangle \log h_k+|\gamma_k|^2h_k^{-1}h_{k+1}-|\gamma_{k-1}|^2h_{k-1}^{-1}h_k&=&0,\quad k=2,\cdots, m-1,\\
\triangle \log h_m+|\mu|^2h_m^{-1}-|\gamma_{m-1}|^2h_{m-1}^{-1}h_m&=&0.
\end{\eq}
The pullback metric is
$g=2n(|\nu|^2h_{1}^{2}+2|\mu|^2h_m^{-1}+2\sum_{k=1}^{m-1}|\gamma_k|^2h_k^{-1}h_{k+1})dz\otimes d\bar{z}.$

\subsection{Hitchin fibration and cyclic Higgs bundles in the Hitchin component}\label{s3}
Fix a Riemann surface $\Sigma$, the Hitchin fibration is a map
\[h:M_{Higgs}(SL(n,\mathbb{C})) \longrightarrow\bigoplus\limits_{j=2}^nH^0(\Sigma, K^j)\ni (q_2,q_3,\cdots,q_n)\]
given by $h([E,\phi])=(\text{tr}(\phi^2),\dots,\text{tr}(\phi^n)).$

In \cite{Hitchin92}, Hitchin defines a section $s_h$ of this fibration whose image consists of stable Higgs bundles with corresponding flat connections having holonomy in $SL(n,\mathbb{R}).$
Furthermore, the section $s_h$ maps surjectively to the connected component (called Hitchin component) of the $SL(n,\mathbb{R})$-Higgs bundle moduli space which naturally contains an embedded copy of Teichm\"uller space. The Teichm\"uller locus is corresponding to the image of $q_3=\cdots=q_n=0$. Such a $(E,\phi)$ corresponds to a representation $\rho$ which can be factored through $SL(2,\mathbb{R})$,
\begin{\eq}
\rho:\pi_1\to SL(2,\mathbb{R})\xrightarrow{\iota} SL(n,\mathbb{R})\hookrightarrow SL(n,\mathbb{C}),
\end{\eq}
where $\iota$ is the canonical irreducible representation.

The cyclic Higgs bundles in the Hitchin component are corresponding to the image of $s_h$ at $(0,\cdots,0,n\cdot q_n)$. More precisely
\begin{\eq}
E=K^{\frac{n-1}{2}}\oplus K^{\frac{n-3}{2}}\oplus\cdots \oplus K^{\frac{3-n}{2}}\oplus K^{\frac{1-n}{2}},\quad
\phi=
\left(
\begin{array}{cccccc}
0 & & & q_n\\
1 & 0 & &\\
& \ddots & \ddots &\\
& & 1 & 0
\end{array}
\right),
\end{\eq}
where $q_n$ is a holomorphic $n$-differential.

If $q_n=0$, the Higgs bundle is Fuchsian. For $n=2m$,
\begin{\eq}
h^{-1}_{k}h_{k+1}=\frac{1}{2}k(n-k)g_{0},\quad 1\leq k\leq m-1,\qquad h^{-2}_{m}=\frac{1}{2}m^2g_0;
\end{\eq}
for $n=2m+1$,
\begin{\eq}
h^{-1}_{k}h_{k+1}=\frac{1}{2}k(n-k)g_{0},\quad 1\leq k\leq m-1,\qquad h^{-1}_m=\frac{1}{2}m(m+1)g_0.
\end{\eq}
Here $g_0$ is the hyperbolic metric such that $\triangle \log g_0=g_0$.

\section{Maximum principle for system}\label{MaxP}
The main tool we use in this paper is the following maximum principle for system. We abuse the same notation $g$ to denote both the metric $g(z)dz\otimes d\bar z$ and the local function $g(z)$ on the surface. Define $\triangle_{g}=g^{-1}\triangle$, which is globally defined, called the Laplacian with respect to the metric $gdz\otimes d\bar z$.
\begin{lem}\label{mp}
Let $(\Sigma,g)$ be a closed Riemannian manifold. For each $1\leq i\leq n$, let $u_i$ be a $C^2$ function on $\Sigma\setminus P_i$, where $P_i$ is an isolated subset of $\Sigma$ ($P_i$ can be empty). Suppose $u_i$ approaches to $+\infty$ around $P_i$. Let $P=\bigcup_{i=1}^{n} P_i$. Let $c_{ij}$ be continuous and bounded functions on $\Sigma\setminus P$, $1\leq i,j\leq n$.
Suppose $c_{ij}$ satisfy the following assumptions: in $\Sigma\setminus P$,\\
$(a)$ cooperative: $c_{ij}\geq 0,~ i\neq j$,\\
$(b)$ column diagonally dominant: $\sum_{i=1}^{n}c_{ij}\leq 0,~ 1\leq j\leq n$,\\
$(c)$ fully coupled: the index set $\{1,\cdots,n\}$ cannot be split up in two disjoint nonempty sets $\alpha,\beta$ such that $c_{ij}\equiv 0$ for $i\in\alpha,j\in \beta.$\\
Let $f_i$ be non-positive continuous functions on $\Sigma\setminus P$, $1\leq i\leq n$ and $X$ be a continuous vector field on $\Sigma$. Suppose $u_i$ satisfies
\begin{eqnarray*}
\triangle_{g} u_i+<X,\nabla u_i>+\sum_{j=1}^{n}c_{ij}u_j=f_i \text{ in } \Sigma\setminus P, \quad 1\leq i\leq n.
\end{eqnarray*}
Consider the following conditions:\\
Condition (1) ~$(f_1,\cdots,f_n)\neq (0,\dots,0)$, i.e., there exists $i_{0}\in \{1,\cdots, n\}$, $p_0\in \Sigma\setminus P$, such that $f_{i_0}(p_0)\neq 0$;\\
Condition (2) ~$P$ is nonempty;\\
Condition (3) ~$\sum_{i=1}^n u_i\geq 0$.\\
Then either condition (1) or (2) imply $u_{i}> 0$, $1\leq i\leq n$. And condition (3) implies either $u_{i}> 0$, $1\leq i\leq n$ or $u_{i}\equiv 0$, $1\leq i\leq n$.
\end{lem}
\pf
Let $A=\{1,\cdots,n\}$. For $S\subseteq A$, set $u_S=\Sigma_{i\in S}u_i$. If $S=\phi$, set $u_S=0$. Let $P_S=\bigcup_{i\in S}P_i$. Then
\begin{\eq}
\triangle_{g} u_{S}+<X,\nabla u_S>+\sum_{i\in A}\sum_{l\in S}c_{li}u_i\leq 0\text{ in } \Sigma\setminus P.
\end{\eq}
Then
\begin{\eq}
\triangle_{g} u_{S}+<X,\nabla u_S>+\sum_{j\notin S}\sum_{l\in S}c_{lj}u_j+\sum_{k\in S}\sum_{l\in S}c_{lk}u_{k}\leq 0\text{ in } \Sigma\setminus P.
\end{\eq}
Then for $S\neq \phi$,
\begin{\eq}
\triangle_{g} u_{S}+<X,\nabla u_S>+\sum_{j\notin S}\sum_{l\in S}c_{lj}(u_{\{j\}\cup S}-u_{S})+\sum_{k\in S}(\sum_{l\in S}c_{lk})(u_{S}-u_{S\setminus \{k\}})\leq 0\text{ in } \Sigma\setminus P.
\end{\eq}
Set
\begin{\eq}
b_{S}=\min_{\Sigma} u_{S}, \quad\check{b}_{S}=\min_{j\notin S,k\in S}\{b_{\{j\}\cup S},b_{S\setminus \{k\}}\},\quad b=\min_{S\subseteq A} b_{S},
\end{\eq}
Notice that all these constants are finite. By the assumptions $(a)(b)$, in $\Sigma\setminus P$, $c_{lj}\geq 0$ for $l\in S,j\notin S$, and $\sum_{l\in S}c_{lk}\leq 0$ for $k\in S$, then
\begin{\eq}
\triangle_{g} \check{b}_{S}+<X,\nabla \check{b}_S>+\sum_{j\notin S}\sum_{l\in S}c_{lj}(u_{\{j\}\cup S}-\check{b}_{S})+\sum_{k\in S}(\sum_{l\in S}c_{lk})(\check{b}_{S}-u_{S\setminus \{k\}})\geq 0 \text{ in } \Sigma\setminus P.
\end{\eq}
Then
\begin{\eq}
\triangle_{g} (u_{S}-\check{b}_{S})+<X,\nabla (u_{S}-\check{b}_{S})>+(-\sum_{j\notin S}\sum_{l\in S}c_{lj}+\sum_{k\in S}(\sum_{l\in S}c_{lk}))(u_{S}-\check{b}_{S})\leq 0\text{ in } \Sigma\setminus P.
\end{\eq}

Step 1: We show that under condition (1) and (2), $u_{S}\geq\check{b}_{S}$ for any $S\subset A$; under condition (3), $u_{S}\geq\check{b}_{S}$ for $S\subsetneq A$. In particular, $b_{S}\geq \check{b}_{S}$ for $S\subset A$ under condition (1) and (2) and for $S\subsetneq A$ under condition (3).

If not, since $u_{S}-\check{b}_{S}$ approaches to $+\infty$ around $P_S$ and continuous outside $P_S$, $u_{S}-\check{b}_{S}$ must attain a negative minimum in $\Sigma\setminus P_S$. First, we suppose $u_{S}-\check{b}_{S}$ is not a constant. By the assumptions $(a)(b)$, in $\Sigma\setminus P$,
\begin{\eq}
-\sum_{j\notin S}\sum_{l\in S}c_{lj}+\sum_{k\in S}(\sum_{l\in S}c_{lk})\leq 0.
\end{\eq}
Then by the strong maximum principle for the single equation (see \cite{MaximumPrinciple}), the minimal point $p\notin \Sigma\setminus P$. So $p\in P\setminus P_S$.
Since $P$ is isolated, we consider $p_n\in \Sigma\setminus P$, $p_n\rightarrow p$. Then
\begin{\eq}
\limsup_{p_n\rightarrow p}\Big(\triangle_{g} (u_{S}-\check{b}_{S})+<X,\nabla (u_{S}-\check{b}_{S})>+(-\sum_{j\notin S}\sum_{l\in S}c_{lj}+\sum_{k\in S}(\sum_{l\in S}c_{lk}))(u_{S}-\check{b}_{S})\Big)(p_n)\leq 0.
\end{\eq}
By the continuity,
\begin{\eq}
\lim_{p_n\rightarrow p}\Big(\triangle_{g} (u_{S}-\check{b}_{S})+<X,\nabla (u_{S}-\check{b}_{S})>\Big)(p_n)=\Big(\triangle_{g} (u_{S}-\check{b}_{S})+<X,\nabla (u_{S}-\check{b}_{S})>\Big)(p)\geq 0.
\end{\eq}
So
\begin{\eq}
\limsup_{p_n\rightarrow p}\Big((-\sum_{j\notin S}\sum_{l\in S}c_{lj}+\sum_{k\in S}(\sum_{l\in S}c_{lk}))(u_{S}-\check{b}_{S})\Big)(p_n)\leq 0.
\end{\eq}
If there exists a subsequence $p_{n_k}$ such that $(-\sum_{j\notin S}\sum_{l\in S}c_{lj}+\sum_{k\in S}(\sum_{l\in S}c_{lk}))(p_{n_k})$ approaches to a negative number, then
\begin{\eq}
\lim_{p_{n_k}\rightarrow p}\Big((-\sum_{j\notin S}\sum_{l\in S}c_{lj}+\sum_{k\in S}(\sum_{l\in S}c_{lk}))(u_{S}-\check{b}_{S})\Big)(p_{n_k})>0.
\end{\eq}
Contradiction. Since $P_S$ is isolated, we have $-\sum_{j\notin S}\sum_{l\in S}c_{lj}+\sum_{k\in S}(\sum_{l\in S}c_{lk})$ is continuous in $\Sigma\setminus P_S$. Then
\begin{\eq}
\triangle_{g} (u_{S}-\check{b}_{S})+<X,\nabla (u_{S}-\check{b}_{S})>+(-\sum_{j\notin S}\sum_{l\in S}c_{lj}+\sum_{k\in S}(\sum_{l\in S}c_{lk}))(u_{S}-\check{b}_{S})\leq 0\text{ in } \Sigma\setminus P_S,
\end{\eq}
Then by the strong maximum principle for the single equation, $u_{S}-\check{b}_{S}$ cannot achieve its negative minimum in $\Sigma\setminus P_S$ unless it is a constant. Contradiction. Second, if $u_{S}-\check{b}_{S}$ is a negative constant, then by the assumptions $(a)(b)$, in $\Sigma\setminus P$,
\begin{\eq}
-\sum_{j\notin S}\sum_{l\in S}c_{lj}+\sum_{k\in S}(\sum_{l\in S}c_{lk})\equiv 0.
\end{\eq}
Then in $\Sigma\setminus P$, $\sum_{l\in S}c_{lk}\equiv 0$ for $k\in S$. Then by the assumptions $(a)(b)$, $c_{ij}\equiv 0$ in $\Sigma\setminus P$, for $j\in S, i\notin S$, which is a contradiction to the assumption $(c)$ unless $S=A$. If $S=A$, for condition $(2)$, we have $u_S$ cannot be a constant. And for condition $(1)$, $u_{S}-\check{b}_{S}$ is a negative constant implies $\sum_{i\in A}f_i\equiv 0$, which also gives a contraction. So we obtain $u_{S}\geq \check{b}_{S}$ on the whole $\Sigma$. For condition $(3)$, we obtain $u_{S}\geq \check{b}_{S}$ for $S\subsetneq A$. So we finish the claim. \\

Step 2: We show $b=0$.

Since $u_{S}=0$ for $S=\phi$, we have $b \leq 0$. If $b<0$, suppose $b$ is achieved by $S_0$, and $|S_0|$ is the smallest among all minimizers. Then $S_0\neq \phi$. Under condition (1) and (2), $u_{S_0}\geq \check{b}_{S_0}$ is automatically true. Under condition (3), we have $u_A\geq 0$ and hence $S_0\subsetneq A$, $u_{S_0}\geq \check{b}_{S_0}$.

Since $c_{ij}$ are bounded, suppose $-\sum_{j\notin S_0}\sum_{l\in S_0}c_{lj}+\sum_{k\in S_0}(\sum_{l\in S_0}c_{lk})\geq -M$, where $M$ is a positive constant. Then in $\Sigma\setminus P$,
\begin{\eq}
&&\triangle_{g} (u_{S_0}-\check{b}_{S_0})+<X,\nabla (u_{S_0}-\check{b}_{S_0})>-M(u_{S_0}-\check{b}_{S_0})\\
&\leq &-\big(M+(-\sum_{j\notin S_0}\sum_{l\in S_0}c_{lj}+\sum_{k\in S_0}(\sum_{l\in S_0}c_{lk}))\big)(u_{S_0}-\check{b}_{S_0}).
\end{\eq}
We have proved $u_{S_0}-\check{b}_{S_0}\geq 0$. Then by the continuity,
\begin{\eq}
\triangle_{g} (u_{S_0}-\check{b}_{S_0})+<X,\nabla (u_{S_0}-\check{b}_{S_0})>-M(u_{S_0}-\check{b}_{S_0})\leq 0 \text{ in } \Sigma\setminus P_{S_0}.
\end{\eq}
Since $b\leq \check{b}_{S_0}\leq b_{S_0}$ and $u_{S_0}$ achieves $b$, we have $\check{b}_{S_0}=b$. Then by the strong maximum principle, $u_{S_{0}}\equiv\check{b}_{S_0}=b$. Then
\begin{\eq}
\triangle_{g} b+<X,\nabla b>+\sum_{j\notin S_0}\sum_{l\in S_0}c_{lj}(u_{\{j\}\cup S_0}-b)+\sum_{k\in S_0}(\sum_{l\in S_0}c_{lk})(b-u_{S_0\setminus \{k\}})\leq 0\text{ in } \Sigma\setminus P.
\end{\eq}
Then by the assumptions $(a)(b)$,
\begin{\eq}
(\sum_{l\in S_0}c_{lk})(b-u_{S_0\setminus \{k\}})&\equiv&0\text{ in } \Sigma\setminus P, \text{ for } k\in S_0.
\end{\eq}
If $b-u_{S_0\setminus \{k\}}=0$ at one point, then $\check{b}_{S_0\setminus \{k\}}=b$, which is a contradiction since $|S_0|$ is the smallest. So in $\Sigma\setminus P$, $\sum_{l\in S_0}c_{lk}\equiv 0$ for $k\in S_0$. As the argument above, it is a contradiction to the assumption $(c)$. Then we obtain $b=0$, in particular, $u_{i}\geq 0$, $1\leq i\leq n$.\\

Step 3: We finish the proof.

Since $u_{i}\geq 0$, we have in $\Sigma\setminus P$,
\begin{\eq}
\triangle_g u_i+<X,u_i>+c_{ii}u_i\leq 0, \quad 1\leq i\leq n.
\end{\eq}
Then as the argument above, by the strong maximum principle, there exists a subset $Z\subseteq A$, such that $u_i\equiv 0$ for $i\in Z$ and $u_j>0$ for $j\notin Z$. Then for $i\in Z$, in $\Sigma\setminus P$, $0\leq \sum_{j\notin Z}c_{ij}u_j=f_i\leq 0$. Since $u_j>0$ for $j\notin Z$, $c_{ij}\equiv 0$ for $i\in Z,~j\notin Z$. Suppose condition (1) $(f_1,\cdots,f_n)\not\equiv (0,\dots,0)$ or condition (2) $P$ is nonempty holds, we can rule out the possibility $Z=A$. Suppose condition (3) $\sum_{i=1}^n u_i\geq 0$ holds, $Z$ must be empty or $A$. So either $u_{i}>0$, $1\leq i\leq n$ or $u_i\equiv 0$ for $1\leq i\leq n$ .
\qed
\begin{rem}
Let $\lambda_i$ be positive numbers, $i=1,\cdots, n$. Let $u_i'=\lambda_iu_i$, $c_{ij}'=c_{ij}\lambda_i\lambda_j^{-1}$. If $c_{ij}'$ satisfy the assumptions $(a)(b)(c)$, then we still obtain the same results for $u_i$.
\end{rem}
\begin{rem}
The assumption $(c)$ is easy to check by the following procedure. If $1\in\alpha$, consider $\beta_1=\{j: c_{1j}\equiv 0\}$, $\alpha_1=\{1,\cdots,n\}\setminus \beta_1$. Then $\alpha_1\cap \beta=\phi$. Then $\alpha_1\subseteq \alpha$. Denote $\alpha_0=\{1\}$. If $\alpha_{1}\subseteq\alpha_0$, then $\alpha=\alpha_0$ gives such a partition. If $\alpha_{1}\nsubseteq\alpha_0$, consider $\beta_2=\{j: c_{ij}\equiv 0, i\in \alpha_{0}\cup\alpha_1 \}$, $\alpha_2=\{1,\cdots,n\}\setminus \beta_2$. Then $\alpha_2\subseteq \alpha$. If $\alpha_2\subseteq \alpha_0\cup\alpha_1$, then $\alpha=\alpha_0\cup\alpha_1$ gives such a partition. If $\alpha_{2}\nsubseteq\alpha_0\cup\alpha_1$, consider $\beta_3=\{j: c_{ij}\equiv 0, i\in \bigcup_{k=0}^{2}\alpha_k \}$, $\alpha_3=\{1,\cdots,n\}\setminus \beta_3$. Repeat this procedure, then either we obtain a partition $\alpha,\beta$ such that $c_{ij}\equiv 0$ for $i\in\alpha,j\in\beta$ or we show that $1\notin \alpha$. If $1\notin \alpha$, repeat the procedure above for $2,3,\cdots, n$. Then we can show whether such a partition exists or not.
\end{rem}

\begin{rem}
The maximum principle above may be applied to the non-linear version under certain assumptions, by using the linearization
\begin{\eq}
F(u_1,\cdots,u_n,x)-F(v_1,\cdots,v_n,x)=\sum_{j=1}^{n}(u_j-v_j)\int_{0}^{1}\frac{\partial F}{\partial u_j}(tu_1+(1-t)v_1,\cdots,tu_n+(1-t)v_n,x)dt.
\end{\eq}
For the problems involving poles, we need to check whether the coefficient after linearization is bounded.
\end{rem}

\section{Monotonicity of pullback metrics}\label{monotonicity}
In this section, we first consider the family of the cyclic Higgs bundles $(E,\phi^t)$ parametrized by $(\ga_1,\cdots,\ga_{n-1},t\ga_n),$ $n\geq 3$ for $t\in \mathbb{C}$. We show the monotonicity of the pullback metrics of the corresponding branched minimal immersions along the family $\phi^t$.

\begin{prop}\label{domination1article}
Let $(E,\phi^t)$ be a family of cyclic Higgs bundles parametrized by $(\ga_1,\cdots,\ga_{n-1},t\ga_n)$, $n\geq 3$, $\gamma_n\neq 0$, $t\in\mathbb{C}^*$ and $h^t$ be the corresponding harmonic metrics on $E$. Then as $|t|$ increases, $h_k^{-1}h_{k+1}$, $k=1,\cdots,n-1$ and $t^2h_n^{-1}h_1$ strictly increase.
As a result, outside the branch points, the pullback metric $g^t$ of the corresponding branched minimal immersions strictly increases.
\end{prop}
\pf
We show that for $0<|t'|<|t|$, all the terms for $t$ dominate the corresponding terms for $t'$.

Let $u_k=h_k^{-1}h_{k+1}$, $k=1,\cdots, n-1$, $u_n=|t|^2h_n^{-1}h_1$. Then
\begin{\eq}
\triangle \log u_k+|\ga_{k+1}|^2u_{k+1}-2|\ga_k|^2u_k+|\ga_{k-1}|^2u_{k-1}&=&0,\quad k=1,\cdots, n,
\end{\eq}

And $\ti{u}_k$ are similarly defined for $t'$, satisfying
\begin{\eq}
\triangle \log \ti{u}_k+|\ga_{k+1}|^2\ti{u}_{k+1}-2|\ga_k|^2\ti{u}_k+|\ga_{k-1}|^2\ti{u}_{k-1}&=&0,\quad k=1,\cdots, n,
\end{\eq}
Let $v_k=\log(u_k\ti{u}^{-1}_{k})$, then
\begin{\eq}
\triangle v_k+|\ga_{k+1}|^2\ti{u}_{k+1}(e^{v_{k+1}}-1)-2|\ga_k|^2\ti{u}_{k}(e^{v_k}-1)+|\ga_{k-1}|^2\ti{u}_{k-1}(e^{v_{k-1}}-1)&=&0,\quad k=1,\cdots, n,
\end{\eq}
Let $$c_k=g_0^{-1}|\gamma_k|^2\ti{u}_k\int_{0}^{1}e^{(1-t)(v_k)}dt,~ k=1,\cdots, n.$$
Then $v_k$'s satisfy
\begin{\eq}
\triangle_{g_0} v_k+c_{k-1}v_{k-1}-2c_kv_k+c_{k+1}v_{k+1}&=&0,\quad k=1,\cdots,n
\end{\eq}
It is easy to check that the above system of equations satisfies the assumptions in Lemma \ref{mp} and condition (3), since $\sum_{k=1}^nv_k=2\log (\frac{|t|}{|t'|})>0$. One can apply the maximum principle Lemma \ref{mp}, then $v_k>0$, $k=1,\cdots, n$. Then we obtain $u_k>\ti{u}_k, ~k=1,\cdots,n.$

Finally, the monotonicity of $g^t$ follows from $g^t=2n(\sum_{k=1}^{n-1}|\gamma_k|^2h_k^{-1}h_{k+1}+|\ga_n|^2t^2h_n^{-1}h_1)dz\otimes d\bar{z}$.
\qed

For $t\in \mathbb{C}^*$, the family $(E,\phi^t)$ is gauge equivalent to $t^{\frac{1}{n}}(E,\phi)=(E,t^{\frac{1}{n}}\phi)$ by the gauge transformation $\psi_t=\text{diag}(t^{\frac{n-1}{2n}},t^{\frac{n-3}{2n}},\cdots,t^{\frac{3-n}{2n}},t^{\frac{1-n}{2n}}),$
since
\begin{\eq}t^{\frac{1}{n}}\begin{pmatrix}
0 & & & \gamma_n\\
\gamma_1 & 0 & & \\
&\ddots & \ddots & \\
& & \gamma_{n-1} & 0
\end{pmatrix}&=&\psi_t^{-1}\begin{pmatrix}
0 & & & t\gamma_n\\
\gamma_1 & 0 & & \\
&\ddots & \ddots & \\
& & \gamma_{n-1} & 0
\end{pmatrix}\psi_t.\end{\eq}
Then we obtain the following results.
\begin{cor}
Let $(E,\phi)$ be a cyclic Higgs bundle parametrized by $(\ga_1,\cdots,\ga_n)$, $n\geq 3$. Let $g^t$ be the pullback metric corresponding to $t\phi$ for $t\in \mathbb{C}^*$. Then outside the branch points, along the $\mathbb{C}^*$-orbit, $g^t$ strictly increases as $|t|$ increases.
\end{cor}
Consider the Morse function $f$ on the moduli space of Higgs bundles as the $L^2$-norm of $\phi$: \begin{\eq} f(E,\phi)=\int_{\Sigma}\text{tr}(\phi\phi^*) \sqrt{-1}dz\wedge d\bar z.\end{\eq}
\begin{cor}
Let $(E,\phi)$ be a cyclic Higgs bundle. Then along the $\mathbb{C}^*$-orbit of $(E,\phi)$, the Morse function $f(E,t\phi)$ strictly increases as $|t|$ increases.
\end{cor}
Applying Proposition \ref{domination1article} to $SL(n,\mathbb{R})$ case, we obtain the monotonicity of the harmonic metric.
\begin{cor}
Let $(E,\phi)$ be a cyclic $SL(n,\mathbb{R})$ Higgs bundle parameterized by $(\nu,\ga_1,\cdots,\ga_{m-1},\mu)$, $\nu\neq 0$. Denote $\nu=\ga_0$, $\mu=\ga_m$. Consider a family of $SL(n,\mathbb{R})$ cyclic Higgs bundles parameterized by $(\ga_0,\cdots,t\ga_l,\cdots,\ga_{m})$, $l=0,\cdots,m$ for $t\in \mathbb{C}^*$. Let $h^{t\ga_l}$ be the corresponding harmonic metrics. Then as $|t|$ increases, $h_{k}^{t\ga_l}$ strictly increases for $k=1,\cdots,l$ and $h_{k}^{t\ga_l}$ strictly decreases for $k=l+1,\cdots, m$.
\end{cor}
If the cyclic Higgs bundles parametrized by $(\ga_1,\cdots,\ga_{n-1},0)$ is stable, we can extend the monotonicity of the pullback metric of $\mathbb{C}^*$-family to $\mathbb{C}$-family.
\begin{prop}\label{domination2article}
Let $(E,\phi)$ be a family of cyclic Higgs bundles parametrized by $(\ga_1,\cdots,\ga_{n-1},\ga_n)$, $n\geq 3$, $\gamma_n\neq 0$ and $h$ be the corresponding harmonic metrics on $E$. If $(E,\tilde\phi)$ be a family of cyclic Higgs bundles parametrized by $(\ga_1,\cdots,\ga_{n-1},0)$ is stable, then $h_k^{-1}h_{k+1}$, $k=1,\cdots,n-1$ and $h_n^{-1}h_1$ for $(E,\phi)$ strictly dominate the items for $(E,\tilde\phi)$.
As a result, outside the branch points, the pullback metric $g$ of the corresponding branched minimal immersions for $(E,\phi)$ strictly dominates the one for $(E,\tilde\phi)$.
\end{prop}
\begin{pf}
Set $n+1=1$, then the equation for $h_k$ is
\begin{\eq}
\triangle \log h_k+|\gamma_k|^2h_k^{-1}h_{k+1}-|\gamma_{k-1}|^2h_{k-1}^{-1}h_k&=&0,\quad k=1,\cdots, n.
\end{\eq}
Let $u_k=h_k^{-1}h_{k+1}$, $k=1,\cdots, n$. Then
\begin{\eq}
\triangle \log u_1+|\gamma_2|^2u_2-2|\ga_1|^2u_1&=&-|\ga_n|^2u_n \leq 0,\\
\triangle \log u_k+|\ga_{k+1}|^2u_{k+1}-2|\ga_k|^2u_k+|\ga_{k-1}|^2u_{k-1}&=&0,\quad k=2,\cdots, n-2,\\
\triangle \log u_{n-1}-2|\ga_{n-1}|^2u_{n-1}+|\ga_{n-2}|^2u_{n-2}&=&-|\ga_{n}|^2u_{n} \leq 0.
\end{\eq}
And $\ti{h}_k,\ti{u}_k$ are similarly defined for $t=0$.

Let $v_k=\log(u_k\ti{u}^{-1}_{k})$, $ k=1,\cdots, n-1$. Then
\begin{\eq}
\triangle v_1+|\gamma_2|^2\ti{u}_{2}(e^{v_2}-1)-2|\ga_1|^2\ti{u}_1(e^{v_1}-1)&\leq& 0,\\
\triangle v_k+|\ga_{k+1}|^2\ti{u}_{k+1}(e^{v_{k+1}}-1)-2|\ga_k|^2\ti{u}_{k}(e^{v_k}-1)+|\ga_{k-1}|^2\ti{u}_{k-1}(e^{v_{k-1}}-1)&=&0,\quad k=2,\cdots, n-2,\\
\triangle v_{n-1}-2|\ga_{n-1}|^2\ti{u}_{n-1}(e^{v_{n-1}}-1)+|\ga_{n-2}|^2\ti{u}_{n-2}(e^{v_{n-2}}-1)&\leq& 0.
\end{\eq}
Let $c_k=g_0^{-1}|\gamma_k|^2\ti{u}_k\int_{0}^{1}e^{(1-t)v_k}dt,~ k=1,\cdots, n-1.$
Then $v_k$'s satisfy
\begin{\eq}
\triangle_{g_0} v_1-2c_1v_1+c_2v_2&\leq&0,\\
\triangle_{g_0} v_k+c_{k-1}v_{k-1}-2c_kv_k+c_{k+1}v_{k+1}&=&0,\quad k=2,\cdots,n-2\\
\triangle_{g_0} v_{n-1}+c_{n-2}v_{n-2}-2c_{n-1}v_{n-1}&\leq&0.
\end{\eq}
It is easy to check that the above system of equations satisfies the assumptions in Lemma \ref{mp} and condition (1), since $\gamma_n\neq 0$. Applying the maximum principle Lemma \ref{mp}, $v_k>0$, $k=1,\cdots, n-1$. Then we obtain $u_k>\ti{u}_k, ~k=1,\cdots,n-1.$
\end{pf}\qed
\begin{rem}
Proposition \ref{domination1article}, \ref{domination2article} is a generalization of the metric domination theorem in \cite{DaiLi} in two aspects: (1) from dominating the Fuchsian case to monotonicity along the $\mathbb{C}$-family; (2) from cyclic Higgs bundles in the Hitchin component to general cyclic Higgs bundles.
\end{rem}

\section{Curvature of cyclic Higgs bundles in the Hitchin component}\label{curvature}
In this section, we would like to obtain a lower and upper bound for the extrinsic curvature of the branched minimal immersion associated to cyclic Higgs bundles.
Let's first get to know how big the range of the sectional curvature of the symmetric space is.
\begin{prop} \label{CurvatureSym}
Let $G=SL(n,\mathbb{C}), SL(n,\mathbb{R}), Sp(2m,\mathbb{R})(n=2m)$, the maximal compact subgroup $K=SU(n), SO(n), U(m)$ respectively. For any tangent plane $\sigma$ in $G/K$, the sectional curvature $K_{\sigma}$ for the associated symmetric space $G/K$ satisfies $$-\frac{1}{n}\leq K_{\sigma}\leq 0,$$ where
(1) for $SL(n,\mathbb{C}), SL(n,\mathbb{R})$, $-\frac{1}{n}$ can be achieved by the tangent plane spanned by $$E_{ij}+E_{ji}, E_{ii}-E_{jj} \quad\text{ for any $1\leq i<j\leq n.$}$$
(2) for $Sp(2m,\mathbb{R})$ where $n=2m$, $-\frac{1}{n}$ can be achieved by the tangent plane spanned by $$E_{i,m+i}+E_{m+i,i},E_{ii}-E_{m+i,m+i}\quad\text{for any~} 1\leq i\leq m.$$
\end{prop}

\pf Suppose the Cartan decomposition of the Lie algebra is $\mathfrak{g}=\mathfrak{k}+\mathfrak{p}$.
The sectional curvature of the plane spanned by the vectors $Y_1,Y_2\in T_p(G/K)$ is (see \cite{Jost} for reference)
\begin{\eq}
K(Y_1\wedge Y_2)=\frac{B([\omega^{\mathfrak{p}}(Y_1),\omega^{\mathfrak{p}}(Y_2)],[\omega^{\mathfrak{p}}(Y_1),\omega^{\mathfrak{p}}(Y_2)])}{B(\omega^{\mathfrak{p}}(Y_1),\omega^{\mathfrak{p}}(Y_1))\cdot B(\omega^{\mathfrak{p}}(Y_2),\omega^{\mathfrak{p}}(Y_2))-B(\omega^{\mathfrak{p}}(Y_1),\omega^{\mathfrak{p}}(Y_2))^2}.
\end{\eq}
So it is enough by only checking $Y_1,Y_2\in T_{eK}(G/K)=\mathfrak{p}$.
The upper bound is obvious since $B$ is negative definite on $\mathfrak{k}$, where $[\omega^{\mathfrak{p}}(Y_1),\omega^{\mathfrak{p}}(Y_2)]$ lies.

Now we show the lower bound.
Let $\sigma$ be the plane $\sigma=\text{span}\{Y,Z\}$ where $Y,Z\in\mathfrak{p}$ satisfying $\text{tr}(YZ)=0, \text{tr}(Y^2)=\text{tr}(Z^2)$. The Killing form $B(Y,Z)=2n\cdot\text{tr}(YZ)$. Define $U=Y+iZ, V=Y-iZ$, then the sectional curvature of the plane $\sigma=\text{span}\{Y,Z\}$ is
\begin{\eq}
K_{\sigma}&=&-\frac{\text{tr}(UVUV)-\text{tr}(U^2V^2)}{n\cdot \text{tr}(UV)^2}\\
&\geq&-\frac{\text{tr}(UVUV)}{n\cdot \text{tr}(UV)^2}\quad \text{using $\text{tr}(U^2V^2)=\text{tr}(U^2\overline{U^2}^T)\geq 0$}\\
&\geq& -\frac{1}{n}\quad\text{using $\text{tr}(A^2)\leq \text{tr}(A)^2.$}
\end{\eq}

The equality holds if and only if $U^2=0$ and $UV=U\overline{U}^T$ has only one nonzero eigenvalue. In terms of $Y,Z$, the equality holds if and only if
$Y^2=Z^2$, $YZ+ZY=0$, and $Y^2+Z^2+i(ZY-YZ)$ has only one nonzero eigenvalue. The rest is by direct calculation.
\qed\\

For general cyclic Higgs bundles, one should not expect a nontrivial lower bound of the extrinsic curvature at immersed points since it could achieve the plane of the most negative curvature in $SL(n,\mathbb{C})/SU(n)$.
\begin{prop} For cyclic Higgs bundles parametrized by $(\ga_1,\ga_2,\cdots,\ga_n)$, if there exists a point such that $n-1$ terms of $\ga_i$'s are equal to zero, the sectional curvature of the tangent plane of the associated harmonic map at this point is $-\frac{1}{n}$.
\end{prop}
\pf Firstly, $n=2$ case is obvious. Let $n\geq 3$. The associated harmonic map is a possibly branched minimal immersion.
The tangent plane $\sigma$ of the minimal immersion at $f(p)$ inside $G/K$ is spanned by $Y_{f(p)}=f_*(\frac{\partial}{\partial x})$ and $Z_{f(p)}=f_*(\frac{\partial}{\partial y})$. Using the formula (\ref{KeyFormula}) in Section \ref{pre},
\begin{\eq}
&&\omega^{\mathfrak{p}}(Y)=\Phi(\frac{\partial}{\partial z})=(\phi+\phi^*)(\frac{\partial}{\partial x})=\phi(\frac{\partial}{\partial z})+\phi^*(\frac{\partial}{\partial \bar{z}}),\\
&&\omega^{\mathfrak{p}}(Z)=\Phi(\frac{\partial}{\partial y})=(\phi+\phi^*)(\frac{\partial}{\partial y})=\sqrt{-1}\phi(\frac{\partial}{\partial z})-\sqrt{-1}\phi^*(\frac{\partial}{\partial \bar{z}}).
\end{\eq}
One may refer the details in Section 2 in \cite{DaiLi}. Hence \begin{\eq}
[\omega^{\mathfrak{p}}(Y),\omega^{\mathfrak{p}}(Z)]=-2\sqrt{-1}[\phi(\frac{\partial}{\partial z}),\phi^*(\frac{\partial}{\partial\bar z})]=-2\sqrt{-1}[\phi,\phi^*](\frac{\partial}{\partial z}, \frac{\partial}{\partial\bar z})
\end{\eq}
Since $f$ is conformal, we have $Y\perp Z$. Then the sectional curvature of the plane $\sigma$ is
\begin{eqnarray}\label{CurvatureFormula}
K_{\sigma}&=& K(Y\wedge Z)=\frac{B([\omega^{\mathfrak{p}}(Y),\omega^{\mathfrak{p}}(Z)],[\omega^{\mathfrak{p}}(Y),\omega^{\mathfrak{p}}(Z)])}{B(\omega^{\mathfrak{p}}(Y),\omega^{\mathfrak{p}}(Y))B(\omega^{\mathfrak{p}}(Z),\omega^{\mathfrak{p}}(Z))}\nonumber\\
&=&-\frac{B([\phi,\phi^*],[\phi,\phi^*])}{B(\phi,\phi^*)B(\phi,\phi^*)}=-\frac{\text{tr}([\phi,\phi^*][\phi,\phi^*])}{2n\cdot\text{tr}(\phi\phi^*)^2}\\
&=&-\frac{(h_n^{-1}h_1|\ga_n|^2-h_1^{-1}h_2|\ga_1|^2)^2+(h_1^{-1}h_2|\ga_1|^2-h_2^{-1}h_3|\ga_2|^2)^2+\cdots
+(h_{n-1}^{-1}h_n|\ga_{n-1}|^2-h_n^{-1}h_1|\ga_n|^2)^2}{2n(h_n^{-1}h_1|\ga_n|^2+h_1^{-1}h_2|\ga_1|^2
+\cdots+h_{n-1}^{-1}h_n|\ga_{n-1}|^2)^2}.\nonumber
\end{eqnarray}
In particular, if at point $p$, there exists $k_0$ such that $\ga_i=0$, for $i\neq k_0$, and $\ga_{k_0}\neq 0$.
Then $$K_p=-\frac{2(h_{k_0-1}^{-1}h_{k_0}|\ga_{k_0}|^2)^2}{2n\cdot(h_{k_0-1}^{-1}h_{k_0})^2|\ga_{k_0}|^2)^2}
=-\frac{1}{n}.$$
\qed
\begin{rem}
For example, consider the cyclic Higgs bundle $(L\oplus\mathcal{O}\oplus L^{-1}, \begin{pmatrix}0&0&\beta\\\alpha&0&0\\0&\alpha&0\end{pmatrix})$, where $\deg L<\deg K, 0\neq\alpha\in H^0(L^{-1}K), 0\neq\beta\in H^0(L^2K)$. Suppose in addition, zeros of $\beta$ do not contain all zeros of $\alpha$. Then at any point where $\alpha=0,\beta\neq 0$, the map is an immersion and the extrinsic curvature is $-\frac{1}{3}$.
\end{rem}
So instead, we restrict ourselves to cyclic Higgs bundles in the Hitchin components. In this case, we obtain a nontrivial lower and upper bound on the extrinsic curvature of the associated minimal immersion into $G/K$.

Let $(E,\phi)$ be a cyclic Higgs bundle in the Hitchin component parameterized by $q_n\neq 0$ and $(E,\ti{\phi)}$ be the Fuchsian case. Let $h,\ti{h}$ be the corresponding harmonic metrics.
For $n=2m$ even, define \begin{\eq}
\nu_1=\frac{h_1^2|q_n|^2}{h_1^{-1}h_2},\quad \nu_k=\frac{h_{k-1}^{-1}h_k}{h_{k}^{-1}h_{k+1}},\quad k=2,\cdots, m-1,\quad \nu_m=\frac{h_{m-1}^{-1}h_m}{h_m^{-2}}.
\end{\eq}

Similarly, define
\begin{\eq}\ti\nu_1=0,\quad \ti\nu_k=\frac{\ti h_{k-1}^{-1}\ti h_k}{\ti h_k^{-1}\ti h_{k+1}},\quad k=2,\cdots,m-1,\quad\ti\nu_m=\frac{\ti h_{m-1}^{-1}\ti h_m}{\ti h_m^{-2}}.
\end{\eq} By the explicit description of $\ti h$, $\ti h_k^{-1}\ti h_{k+1}=\frac{1}{2}k(n-k)g_0$ for $k=1,\cdots,m-1$ and $\ti h^{-2}_{m}=\frac{m(n-m)}{2}g_0$. Here $g_0$ is the hyperbolic metric.

For $n=2m+1$ odd, $\nu_k,\ti\nu_k, \ti{h}_k$ are as above except $\nu_m=\frac{h_{m-1}^{-1}h_m}{h_m^{-1}}$, $\ti\nu_m=\frac{\ti h_{m-1}^{-1}\ti h_m}{\ti h_m^{-1}}$, $\ti h^{-1}_m=\frac{m(n-m)}{2}g_0$.
\begin{prop}\label{keylemma} In the above settings,
\begin{\eq}
\frac{(k-1)(n-k+1)}{k(n-k)}=\ti{\nu}_k<\nu_k<1, \quad\quad k=1,\cdots, m.
\end{\eq} 
\end{prop}
\begin{rem}
The inequality $\nu_k<1$ recovers Lemma 5.3, $q_n$ case in \cite{DaiLi}. Here we give a new proof using the maximum principle Lemma \ref{mp} directly.
\end{rem}

\pf
We only prove the case for $n=2m$. The proof is similar for $n=2m+1$.


The equation system for $h_k$ is
\begin{\eq}
\triangle \log h_1+h_1^{-1}h_{2}-|q_n|^2h_{1}^{2}&=&0,\\
\triangle \log h_k+h_k^{-1}h_{k+1}-h_{k-1}^{-1}h_k&=&0,\quad k=2,\cdots, m-1,\\
\triangle \log h_m+h_m^{-2}-h_{m-1}^{-1}h_m&=&0.
\end{\eq}
Let
\begin{\eq}
u_0=\log(|q_n|^2h^2_1), \quad u_k=\log(h^{-1}_kh_{k+1}),\quad 1\leq k\leq m-1,\quad u_m=\log (h^{-2}_m).
\end{\eq}
By the holomorphicity, $\triangle\log |q_n|=0$ outside the zeros of $q_n$. Then outside the zeros of $q_n$,
\begin{\eq}
\triangle u_0+2e^{u_1}-2e^{u_0}&=&0,\\
\triangle u_k+e^{u_{k+1}}-2e^{u_k}+e^{u_{k-1}}&=&0, \quad k=1,\cdots, m-1,\\
\triangle u_m-2e^{u_m}+2e^{u_{m-1}}&=&0.
\end{\eq}

To prove $\nu_k<1$, let $v_k=u_{k+1}-u_k$, $c_k=\int_0^1e^{tu_{k+1}+(1-t)u_k}dt$, $k=0,\cdots, m-1$. Then outside the zeros of $q_n$,
\begin{\eq}
\triangle v_0-3c_0v_0+c_1v_1&=&0, \\
\triangle v_k+c_{k-1}v_{k-1}-2c_kv_k+c_{k+1}v_{k+1}&=&0, \quad k=1,\cdots, m-2,\\
\triangle v_{m-1}+c_{m-2}v_{m-2}-3c_{m-1}v_{m-1}&=&0.
\end{\eq}
Note that only $v_0$ has poles at zeros of $q_n$. To apply Lemma \ref{mp}, we check that $c_0$ is bounded. In fact, around the zeros of $q_n$,
\begin{\eq}
c_0=\int_{0}^{1}e^{tu_1+(1-t)u_0}dt=\int_{0}^{1}(|q_n|^2h_1^2)^{1-t}e^{tu_1}dt
\leq C.
\end{\eq}
It is then easy to check that the above system of equations satisfies the assumptions in Lemma \ref{mp} and condition (2), since the set of poles (i.e. the set of zeros of $q_n$) is nonempty. Applying the maximum principle Lemma \ref{mp}, $v_k>0$, $k=1,\cdots, m-1$. Then we obtain $\nu_k<1, ~k=1,\cdots, m-1.$\\

To prove $\nu_k>\ti{\nu}_k$, define
\begin{\eq}
u_k=\log(h^{-1}_kh_{k+1}), \quad 1\leq k\leq m-1, \quad u_m=\log (h^{-2}_m).
\end{\eq}
Then
\begin{\eq}
\triangle (u_{2}-u_1)+e^{u_{3}}-3e^{u_{2}}+3e^{u_{1}}&=&|q_n|^2h_1^2\geq0,\\
\triangle (u_{k+1}-u_k)+e^{u_{k+2}}-3e^{u_{k+1}}+3e^{u_{k}}-e^{u_{k-1}}&=&0,\quad k=2,\cdots, m-2,\\
\triangle (u_{m}-u_{m-1})-3e^{u_{m}}+4e^{u_{m-1}}-e^{u_{m-2}}&=&0.
\end{\eq}
And $\ti{u}_k$ are similarly defined for the Fuchsian case, satisfying
\begin{\eq}
\triangle (\ti u_{2}-\ti u_1)+e^{\ti u_{3}}-3e^{\ti u_{2}}+3e^{\ti u_{1}}&=& 0,\\
\triangle (\ti u_{k+1}-\ti u_k)+e^{\ti u_{k+2}}-3e^{\ti u_{k+1}}+3e^{\ti u_{k}}-e^{\ti u_{k-1}}&=&0,\quad k=2,\cdots, m-2,\\
\triangle (\ti u_{m}-\ti u_{m-1})-3e^{\ti u_{m}}+4e^{\ti u_{m-1}}-e^{\ti u_{m-2}}&=&0.
\end{\eq}
To estimate $(u_{k+1}-u_k)-(\ti{u}_{k+1}-\ti{u}_k)$, we have for $k=2,\cdots, m-2$,
\begin{\eq}
&&(e^{u_{k+2}}-3e^{u_{k+1}}+3e^{u_{k}}-e^{u_{k-1}})-(e^{\ti{u}_{k+2}}-3e^{\ti{u}_{k+1}}+3e^{\ti{u}_{k}}-e^{\ti{u}_{k-1}})\\
&=&e^{\ti{u}_{k+2}}(e^{u_{k+2}-\ti{u}_{k+2}}-e^{u_{k+1}-\ti{u}_{k+1}})-2e^{\ti{u}_{k+1}}(e^{u_{k+1}-\ti{u}_{k+1}}-e^{u_{k}-\ti{u}_{k}})
+e^{\ti{u}_{k}}(e^{u_{k}-\ti{u}_{k}}-e^{u_{k-1}-\ti{u}_{k-1}})\\
&&+(e^{\ti{u}_{k+2}}-e^{\ti{u}_{k+1}})(e^{u_{k+1}-\ti{u}_{k+1}}-1)-2(e^{\ti{u}_{k+1}}-e^{\ti{u}_{k}})(e^{u_{k}-\ti{u}_{k}}-1)
+(e^{\ti{u}_{k}}-e^{\ti{u}_{k-1}})(e^{u_{k-1}-\ti{u}_{k-1}}-1)\\
&=&e^{\ti{u}_{k+2}}(e^{u_{k+2}-\ti{u}_{k+2}}-e^{u_{k+1}-\ti{u}_{k+1}})-2e^{\ti{u}_{k+1}}(e^{u_{k+1}-\ti{u}_{k+1}}-e^{u_{k}-\ti{u}_{k}})
+e^{\ti{u}_{k}}(e^{u_{k}-\ti{u}_{k}}-e^{u_{k-1}-\ti{u}_{k-1}})\\
&&+(e^{u_{k+1}-\ti{u}_{k+1}}-e^{u_{k}-\ti{u}_{k}})(e^{\ti{u}_{k+2}}-e^{\ti{u}_{k+1}})-(e^{u_{k}-\ti{u}_{k}}-e^{u_{k-1}-\ti{u}_{k-1}})(e^{\ti{u}_{k}}-e^{\ti{u}_{k-1}})\\
&&+(e^{u_k-\ti{u}_{k}}-1)(e^{\ti{u}_{k+2}}-3e^{\ti{u}_{k+1}}+3e^{\ti{u}_{k}}-e^{\ti{u}_{k-1}}).
\end{\eq}
Since $\ti{u}_{k+1}-\ti{u}_k$ is a globally defined constant function, the equation of $\ti{u}_{k+1}-\ti{u}_k$ gives
\begin{\eq}
e^{\ti{u}_{k+2}}-3e^{\ti{u}_{k+1}}+3e^{\ti{u}_{k}}-e^{\ti{u}_{k-1}}=0.
\end{\eq}
Then
\begin{\eq}
&&(e^{u_{k+2}}-3e^{u_{k+1}}+3e^{u_{k}}-e^{u_{k-1}})-(e^{\ti{u}_{k+2}}-3e^{\ti{u}_{k+1}}+3e^{\ti{u}_{k}}-e^{\ti{u}_{k-1}})\\
&=&e^{\ti{u}_{k+2}}(e^{u_{k+2}-\ti{u}_{k+2}}-e^{u_{k+1}-\ti{u}_{k+1}})+(e^{\ti{u}_{k+2}}-3e^{\ti{u}_{k+1}})(e^{u_{k+1}-\ti{u}_{k+1}}-e^{u_{k}-\ti{u}_{k}})
+e^{\ti{u}_{k-1}}(e^{u_{k}-\ti{u}_{k}}-e^{u_{k-1}-\ti{u}_{k-1}}).
\end{\eq}
Similarly, for $k=1$,
\begin{\eq}
(e^{u_{3}}-3e^{u_{2}}+3e^{u_{1}})-(e^{\ti{u}_{3}}-3e^{\ti{u}_{2}}+3e^{\ti{u}_{1}})
=e^{\ti{u}_{3}}(e^{u_{3}-\ti{u}_{3}}-e^{u_{2}-\ti{u}_{2}})+(e^{\ti{u}_{3}}-3e^{\ti{u}_{2}})(e^{u_{2}-\ti{u}_{2}}-e^{u_{1}-\ti{u}_{1}}),
\end{\eq}
for $k=m-1$,
\begin{\eq}
&&(-3e^{u_{m}}+4e^{u_{m-1}}-e^{u_{m-2}})-(-3e^{\ti{u}_{m}}+4e^{\ti{u}_{m-1}}-e^{\ti{u}_{m-2}})\\
&=&-3e^{\ti{u}_{m}}(e^{u_{m}-\ti{u}_{m}}-e^{u_{m-1}-\ti{u}_{m-1}})+(4e^{\ti{u}_{m-1}}-3e^{\ti{u}_{m}})(e^{u_{m-1}-\ti{u}_{m-1}}-e^{u_{m-2}-\ti{u}_{m-2}}).
\end{\eq}
Let $v_k=(u_{k+1}-u_k)-(\ti{u}_{k+1}-\ti{u}_k)$, $k=1,\cdots, m-1$. Let $c_k=\int_{0}^{1}e^{t(u_{k+1}-\ti{u}_{k+1})+(1-t)(u_{k}-\ti{u}_{k})}dt$, $k=1,\cdots, m-1$. Then
\begin{\eq}
\triangle v_1+e^{\ti{u}_{3}}c_{2}v_{2}+(e^{\ti{u}_{3}}-3e^{\ti{u}_{2}})c_1v_1&\geq&0\\
\triangle v_k+e^{\ti{u}_{k+2}}c_{k+1}v_{k+1}+(e^{\ti{u}_{k+2}}-3e^{\ti{u}_{k+1}})c_kv_k+e^{\ti{u}_{k-1}}c_{k-1}v_{k-1}&=&0,\quad k=2,\cdots, m-2,\\
\triangle v_{m-1}-3e^{\ti{u}_{m}}c_{m-1}v_{m-1}+(4e^{\ti{u}_{m-1}}-3e^{\ti{u}_{m}})c_{m-2}v_{m-2}&=&0.
\end{\eq}
To apply the maximum principle, we need to check
\begin{\eq}
e^{\ti{u}_{k+2}}-2e^{\ti{u}_{k+1}}+e^{\ti{u}_{k}}\leq 0,\quad k=1,\cdots, m-2.
\end{\eq}
This is from the equation of $\ti{u}_{k+1}$ and the fact $\ti{u}_{k+1}=\text{const}+\log g_0$, $\triangle \log g_0=g_0$. Other conditions to apply the maximum principle hold clearly (for $e^{\ti{u}_{m-1}}\leq e^{\ti{u}_m}$, it is from Lemma \ref{keylemma}), so we obtain the desired result.
\qed\\

The cyclic Higgs bundles in the Hitchin component for $n\geq 3$ induce minimal immersions $f:\widetilde{\Sigma}\rightarrow SL(n,\mathbb{R})/SU(n)$. We want to investigate that, as an immersed submanifold, how $f(\widetilde{\Sigma})$ sits in the symmetric space.

\begin{thm}\label{curvaturearticle}
Let $f:\widetilde{\Sigma}\rightarrow SL(n,\mathbb{R})/SU(n)$ be the harmonic map associated to Higgs bundles in the Hitchin component parameterized by $q_n$. Then the sectional curvature $K_{\sigma}$ of the tangent plane $\sigma$ of the image of $f$ in $G/K$ satisfies $$-\frac{1}{n(n-1)^2}\leq K_{\sigma}<0.$$
The equality can be achieved only if $n=2,3$.
\end{thm}

\pf In the case $n=2$, the extrinsic curvature is constantly $-\frac{1}{2}$. Now we consider $n\geq 3$ case. We only prove the case for $n=2m$. The proof is similar for $n=2m+1$.
Using the curvature formula (\ref{CurvatureFormula}), the sectional curvature of the plane $\sigma$ is
\begin{\eq}
K_{\sigma}&=&-\frac{(h_1^2|q_n|^2-h_1^{-1}h_2)^2+(h_1^{-1}h_2-h_2^{-1}h_3)^2+\cdots
+(h_{m-1}^{-1}h_m-h_m^{-2})^2}{n(h_1^2|q_n|^2+2h_1^{-1}h_2+2h_2^{-1}h_3
+\cdots+2h_{m-1}^{-1}h_m+h_m^{-2})^2}.
\end{\eq}
Then $K_\sigma<0$ follows from Proposition \ref{keylemma}.

To show $K_{\sigma}\geq-\frac{1}{n(n-1)^2}$, let $\mu_k=\nu_k^{-1}$, then
\begin{\eq}
K_\sigma &\geq&-\frac{(h_1^{-1}h_2)^2+(h_1^{-1}h_2-h_2^{-1}h_3)^2+\cdots
+(h_{m-1}^{-1}h_m-h_m^{-2})^2}{n(2h_1^{-1}h_2+2h_2^{-1}h_3+\cdots+2h_{m-1}^{-1}h_m+h_m^{-2})^2}\\
&=&-\frac{1+(1-\mu_2)^2+(1-\mu_3)^2\mu_2^{2}+\cdots
+(1-\mu_m)^2\mu_2^{2}\cdots\mu_{m-1}^{2}}{n(2+2\mu_2+2\mu_2\mu_3+\cdots+\mu_2\cdots\mu_{m})^2}
\end{\eq}

Define the functions $G_k, H_k$ for $3\leq k\leq m+1$ as follows. For $3\leq k\leq m-1$,\begin{\eq}
&&G_k=(1-\mu_k)^2+(1-\mu_{k+1})^2\mu_k^2+\cdots+(1-\mu_m)^2\mu_k^2\cdots\mu_{m-1}^2\\
&&H_k=2+2\mu_k+2\mu_k\mu_{k+1}+\cdots+2\mu_k\cdots\mu_{m-1}+\mu_k\cdots\mu_m
\end{\eq}
and \begin{\eq}
G_m=(1-\mu_m)^2,\quad H_m=2+\mu_m,\quad G_{m+1}=0, \quad H_{m+1}=1.
\end{\eq}

The derivatives in $\mu_k$ for $3\leq k\leq m$ are,
\begin{\eq}
(G_k)_{\mu_k}=2\mu_k(1+G_{k+1})-2,\quad (H_k)_{\mu_k}=H_{k+1}
\end{\eq}

Define $F_k$ as a function of $\mu_k$, for $3\leq k\leq m+1$,
\begin{\eq}F_k(\mu_k)=\frac{1+G_k}{(2(k-2)+H_k)^2}.
\end{\eq}
So $K_{\sigma}\geq -\frac{1}{n}F_2.$ For $3\leq k\leq m$,
\begin{\eq}
F_k(\mu_k)=\frac{1+G_k}{(2(k-2)+H_k)^2}=\frac{1+(1-\mu_k)^2+\mu_k^2G_{k+1}}{(2(k-1)
+\mu_kH_{k+1})^2}.
\end{\eq}
We claim:
\begin{lem}\label{lemma1}
$F_2<F_3$.
\end{lem}
\begin{lem}\label{lemma2}
$F_k<F_{k+1}$, for $3\leq k\leq m$.
\end{lem}

Therefore, combining Lemma \ref{lemma1} and \ref{lemma2}, the sectional curvature
\begin{\eq}
K_{\sigma}\geq-\frac{1}{n}F_2>-\frac{1}{n}F_{m+1}=\frac{-1}{n(n-1)^2},\quad\text{for $m>1$}.
\end{\eq}
\qed

\pf (of Lemma \ref{lemma1})
The derivative of $F_2$ in $\mu_2$ is
\begin{\eq}
(F_2)_{\mu_2}&=&\frac{2(2(\mu_2(1+G_{3})-1)-(2-\mu_2)H_{3})}{(1+(1-\mu_2)^2+\mu_2^2G_3)(2+\mu_2H_3)^3}
=\frac{2F}{(1+(1-\mu_2)^2+\mu_2^2G_3)(2+\mu_2H_3)^3},
\end{\eq} where $F=2(\mu_2(1+G_{3})-1)-(2-\mu_2)H_{3}$.
Then
\begin{\eq}
F&<& 2(\ti\mu_2(1+G_{3})-1)-(2-\ti\mu_2)H_{3}\\
&=&2(\frac{2(n-2)}{n-1}(1+G_{3})-1)-(2-\frac{2(n-2)}{n-1})H_{3}\\
&=&\frac{2}{n-1}(n-3+2(n-2)G_{3}-H_{3})\\
&=&\frac{2}{n-1}(n-3+2(n-2)((1-\mu_3)^2+\cdots+(1-\mu_m)^2\mu_3^2\cdots\mu_{m-1}^2)
-(2+2\mu_3+\cdots+2\mu_3\cdots\mu_{m-1})-\mu_3\cdots\mu_m)\\
&=&\frac{2}{n-1}(n-3+P_m),
\end{\eq}
where $P_k=2(n-2)((1-\mu_3)^2+\cdots+(1-\mu_k)^2\mu_3^2\cdots\mu_{k-1}^2)
-(2+2\mu_3+\cdots+2\mu_3\cdots\mu_{k-1})-(n+1-2k)\mu_3\cdots\mu_k$, for $3\leq k\leq m$.

Claim: $P_{k+1}< P_{k}$, for $3\leq k\leq m-1$.
\begin{\eq}
P_{k+1}&=&2(n-2)((1-\mu_3)^2+\cdots+(1-\mu_{k+1})^2\mu_3^2\cdots\mu_{k}^2)
-(2+2\mu_3+\cdots+2\mu_3\cdots\mu_{k})-(n-1-2k)\mu_3\cdots\mu_{k+1}\\
&=&2(n-2)((1-\mu_3)^2+\cdots(1-\mu_{k})^2\mu_3^2\cdots\mu_{k-1}^2)-(2+2\mu_3+\cdots+2\mu_3\cdots\mu_{k})\\
&&+2(n-2)(1-\mu_{k+1})^2\mu_3^2\cdots\mu_{k}^2-(n-1-2k)\mu_3\cdots\mu_{k+1}
\end{\eq}
The last term $2(n-2)(1-\mu_{k+1})^2\mu_3^2\cdots\mu_{k}^2-(n-1-2k)\mu_3\cdots\mu_{k+1}$ satisfies
\begin{\eq}
&&2(n-2)(1-\mu_{k+1})^2\mu_3^2\cdots\mu_{k}^2-(n-1-2k)\mu_3\cdots\mu_{k+1}\nonumber\\
&=&\mu_3\cdots\mu_{k}(2(n-2)\mu_3\cdots\mu_{k}(1-\mu_{k+1})^2-(n-1-2k)\mu_{k+1})\nonumber\\
&<&\mu_3\cdots\mu_{k}(2(n-2)\ti\mu_3\cdots\ti\mu_k(1-\mu_{k+1})^2-(n-1-2k)\mu_{k+1})\nonumber\\
&=&\mu_3\cdots\mu_{k}(k(n-k)(1-\mu_{k+1})^2-(n-1-2k)\mu_{k+1})\nonumber\\
&=&\mu_3\cdots\mu_{k} k(n-k)(\mu_{k+1}^2-(2+\frac{n-1-2k}{k(n-k)})\mu_{k+1}+1)\\
&&\text{Since $1<\mu_{k+1}< \ti\mu_{k+1}=\frac{(k+1)(n-1-k)}{k(n-k)}=1+\frac{n-1-2k}{k(n-k)}$, by Proposition \ref{keylemma}.}\\
&<&-(n-1-2k)\mu_3\cdots\mu_k.
\end{\eq}

Hence $P_{k+1}< P_k$. So \begin{\eq}P_m&<& P_3=2(n-2)(1-\mu_3)^2-2-(n-5)\mu_3\\
&\leq& 2(n-2)(\mu_3^2-(2+\frac{n-5}{2(n-2)})\mu_3+1)< -(n-3).\end{\eq}

Hence $F<0$ and then $(F_2)_{\mu_2}<0$.
Therefore $F_2(\mu_2)< F_2(1)=F_3$.\qed

\pf (of Lemma \ref{lemma2}) The derivative of $F_k$ with respect to $\mu_k$ is
\begin{\eq}
(F_k)_{\mu_k}&=&\frac{2(2(k-1)(\mu_k+\mu_kG_{k+1}-1)
-(2-\mu_k)H_{k+1})}{(2(k-1)+\mu_kH_{k+1})(1+(1-\mu_k)^2+\mu_k^2G_{k+1})^3}
\end{\eq}
By Proposition \ref{keylemma}, $\mu_k< \ti\mu_k=\frac{k(n-k)}{(k-1)(n+1-k)}$,
\begin{\eq}
G_k&<&(1-\ti\mu_k)^2+(1-\ti\mu_{k+1})^2\ti\mu_k^2+\cdots+(1-\ti\mu_m)^2\ti\mu_k^2\cdots\ti\mu_{m-1}^2\\
&=&\frac{1^2+3^2+\cdots+(n+1-2k)^2}{(k-1)^2(n+1-k)^2}=\frac{(n+1-2k)(n+2-2k)(n+3-2k)}{6(k-1)^2(n+1-k)^2}.
\end{\eq}
By Proposition \ref{keylemma}, $\mu_k>1$, then $H_k> n+3-2k.$

The term $2(k-1)(\mu_k-1+\mu_kG_{k+1})-(2-\mu_k)H_{k+1}$ satisfies
\begin{\eq}
&&2(k-1)(\mu_k-1+\mu_kG_{k+1})-(2-\mu_k)H_{k+1}\\
&<& 2(k-1)((\ti\mu_k-1)+\ti\mu_kG_{k+1})-(2-\ti\mu_k)H_{k+1}\\
&<&\frac{2(k-1)}{(k-1)(n+1-k)}((n+1-2k)+\frac{(n-1-2k)(n-2k)(n+1-2k)}{6k(n-k)})\\
&&-(2-\frac{k(n-k)}{(k-1)(n+1-k)})(n+1-2k)\\
&= &\frac{n+1-2k}{(k-1)(n+1-k)}(\frac{(n-1-2k)(n-2k)(k-1)}{3k(n-k)}-(k-2)(n-k))\\
&= &\frac{(n+1-2k)(n-k)}{(k-1)(n+1-k)}(\frac{(n-1-2k)(n-2k)(k-1)}{3(n-k)^2k}-(k-2))\\
&<&\frac{(n+1-2k)(n-k)}{(k-1)(n+1-k)}(\frac{1}{3}-(k-2))<0,\quad\text{for $k\geq 3$.}
\end{\eq}

Hence $F_k$ decreases as $\mu_k$ increases. Then $F_k(\mu_k)< F_k(1)=F_{k+1}$.\qed
\begin{rem}
As shown in \cite{CL14}, along the family of Higgs bundles parameterized by $tq_n$ ($q_n\neq 0$) for $t\in \mathbb{C}$, as $|t|\rightarrow+\infty$, away from zeros of $q_n$, the curvature $K_{\sigma}^t$ approaches to $0$.
\end{rem}
\begin{rem}
The sectional curvature of $SL(n,\mathbb{R})/SO(n)$ satisfies $-\frac{1}{n}\leq K\leq 0$. So the lower bound $-\frac{1}{n(n-1)^2}$ is nontrivial.
\end{rem}
\begin{rem}
(1) In the Fuchsian case, i.e. $q_n=0$, the sectional curvature $K_{\sigma}$ is $-\frac{6}{n^2(n^2-1)}$. Note that $-\frac{6}{n^2(n^2-1)}\geq-\frac{1}{n(n-1)^2}$ and equality holds for $n=2,3$.
\\
(2) At the zeros $p$ of $q_n$, $K_{\sigma}\leq -\frac{6}{n^2(n^2-1)}$ and equality holds if and only if $n=2,3$. For example, in the case $n=2m\geq 4$,
\begin{\eq}
K_{\sigma}(p)&=&-\frac{(h_1^{-1}h_2)^2+(h_1^{-1}h_2-h_2^{-1}h_3)^2+\cdots
+(h_{m-1}^{-1}h_m-h_m^{-2})^2}{n(2h_1^{-1}h_2+2h_2^{-1}h_3+\cdots+2h_{m-1}^{-1}h_m+h_m^{-2})^2}\\
&=&-\frac{(h_1^{-1}h_2)^2+(h_1^{-1}h_2-h_2^{-1}h_3)^2+\cdots
+(h_{m-1}^{-1}h_m-h_m^{-2})^2}{n((2m-1)(h_1^{-1}h_2)
-(2m-3)(h_1^{-1}h_2-h_2^{-1}h_3)-\cdots-(h_{m-1}^{-1}h_m-h_m^{-2}))^2}\\
&&\text{by using the Cauchy--Schwarz inequality and $\nu_k>\ti\nu_k$ for $k=2,\cdots,m$}\\
&< &-\frac{1}{n((2m-1)^2+(2m-3)^2+\cdots+1^2)}\\
&=&-\frac{6}{n^2(n^2-1)}.
\end{\eq}
The case $n=2m+1$ is similar.
\end{rem}

\section{Comparison inside the real Hitchin fibers at $(0,\cdots,0,q_n)$}\label{comparison}
Fix a Riemann surface $\Sigma$, the Hitchin fibration is a map from moduli space of Higgs bundles to the direct sum of holomorphic differentials. We restrict to the $SL(n,\mathbb{R})$-Higgs bundles.

We first compare the harmonic metrics for 
cyclic $SL(n,\mathbb{R})$-Higgs bundles $(E,\phi)$ in the Hitchin fiber at $(0,\cdots,0,n\cdot q_n)$, that is, $\det\phi=(-1)^{n-1}q_n$.
\begin{prop}\label{harmonicmetriccomparison}
Let $(\tilde{E},\tilde{\phi})$ be a cyclic Higgs bundle in the Hitchin component parameterized by $q_n$ and $(E,\phi)$ be a distinct cyclic $SL(n,\mathbb{R})$-Higgs bundle in Section \ref{s2} satisfying $\det\phi=(-1)^{n-1}q_n$. Let $h,\ti{h}$ be the corresponding harmonic metrics.

(1) For $n=2m$, suppose $\ga_1^2\ga_2^2\cdots\ga_{m-1}^2\mu\nu=q_n$, then
\begin{\eq}
h_k>|\ga_k|^2\cdots|\ga_{m-1}|^2|\mu|\ti{h}_k,~k=1,\cdots, m-1,&&\qquad h_m>|\mu|\ti{h}_m.\\
h^{-1}_{m+1-k}>|\nu||\ga_1|^2\cdots|\ga_{m-k}|^2\ti{h}_k,~k=1,\cdots, m-1,&&\qquad h^{-1}_{1}>|\nu|\ti{h}_m.
\end{\eq}

(2) For $n=2m+1$, suppose $\ga_1^2\ga_2^2\cdots\ga_{m-1}^2\mu^2\nu=q_n$, then
\begin{\eq}
h_k>|\ga_k|^2\cdots|\ga_{m-1}|^2|\mu|^2\ti{h}_k,~k=1,\cdots, m-1,&&\qquad h_m>|\mu|^2\ti{h}_m.
\end{\eq}
\end{prop}
\pf
We only prove the inequalities on the first line for $n=2m$. For other cases, the proofs are similar. Define a new Hermitian metric on each $L_k$,
\begin{\eq}
\hat{h}_k=|\ga_k|^2\cdots|\ga_{m-1}|^2|\mu|\ti{h}_k,~k=1,\cdots, m-1, \quad \hat{h}_m=|\mu|\ti{h}_m.
\end{\eq}
By the holomorphicity, $\triangle\log |\ga_k|=0$ outside the zeros of $\ga_k$ (similar for $\mu,\nu$).
Then $\hat{h}$ satisfies, outside the zeros of $q_n$, locally
\begin{\eq}
\triangle \log \hat{h}_1+|\gamma_1|^2\hat{h}_1^{-1}\hat{h}_{2}-|\nu|^2\hat{h}_{1}^{2}&=&0,\\
\triangle \log \hat{h}_k+|\gamma_k|^2\hat{h}_k^{-1}\hat{h}_{k+1}-|\gamma_{k-1}|^2\hat{h}_{k-1}^{-1}\hat{h}_k&=&0, \quad k=2,\cdots, m-1,\\
\triangle \log \hat{h}_m+|\mu|^2\hat{h}_m^{-2}-|\gamma_{m-1}|^2\hat{h}_{m-1}^{-1}\hat{h}_m&=&0.
\end{\eq}
Notice that $\hat{h}$ satisfies the same equation system as $h$, but have zeros.

Define $u_i=\log (h_i/\hat{h}_i)$ and $u_i$ goes to $+\infty$ around the set $P_i$, the zeros of $\hat{h}_i$. Let
\begin{\eq}
c_1&=&g_0^{-1}|\nu|^2\hat{h}_{1}^{2}\int_{0}^{1}e^{(1-t)u_1}dt,\\
c_k&=&g_0^{-1}|\gamma_k|^2\hat{h}_{k-1}^{-1}\hat{h}_{k}\int_{0}^{1}e^{(1-t)u_k}dt,\quad k=2,\cdots, m\\
c_{m+1}&=&g_0^{-1}|\mu|^2\hat{h}_{m}^{-2}\int_{0}^{1}e^{(1-t)u_m}dt.
\end{\eq}
Then $u_i$'s satisfy
\begin{\eq}
\triangle_{g_0} u_1-(c_2+2c_1)u_1+c_2u_2&=&0,\\
\triangle_{g_0} u_k+c_{k+1}u_{k+1}-(c_k+c_{k+1})u_k+c_{k}u_{k-1}&=&0, \quad k=2,\cdots, m-1,\\
\triangle_{g_0} u_m-(2c_{m+1}+c_{m})u_m+c_{m}u_{m-1}&=&0.
\end{\eq}
We need to check the coefficients are bounded. The $c_i$'s are indeed bounded from the fact $\int_{0}^{1}x^{1-t}dt\leq C$ around $x=0$. 
It is then easy to check that the above system of equations satisfies the assumptions in Lemma \ref{mp} and condition (2), since the set $P=\bigcup_i P_i$ of poles is nonempty. Applying Lemma \ref{mp} (the maximum principle), we obtain $u_k>0$, $k=1,\cdots, m$.
\qed\\

Concerning the associated harmonic maps $f:\widetilde{\Sigma}\rightarrow G/K$. We show that the pullback metric of the harmonic map for the cyclic Higgs bundle in the Hitchin component parameterized by $q_n$ dominates the ones for other cyclic $SL(n,\mathbb{R})$-Higgs bundles in the Hitchin fiber at $(0,\cdots, 0, n\cdot q_n)$ for $n=2,3,4$.
\begin{thm}\label{comparisonarticle}
Let $(\tilde{E},\tilde{\phi})$ be a cyclic Higgs bundle in the Hitchin component parameterized by $q_n$ and $(E,\phi)$ be a distinct cyclic $SL(n,\mathbb{R})$-Higgs bundle in Section \ref{s2} such that $\det\phi=(-1)^{n-1}q_n$.
\\
In the case (1) $n=2,\mu\nu=q_2$, (2) $n=3,\mu^2\nu=q_3$, (3) $n=4,\ga^2\mu\nu=q_4$, the pullback metrics $g,\ti{g}$ of corresponding harmonic maps satisfy $g<\ti{g}.$
\end{thm}
\pf
For $n=2$, locally
\begin{\eq}
\frac{1}{4}g=q_2dz^2+(|\nu|^2h^2+|\mu|^2h^{-2})dz\cdot d\bar z+\bar q_2d\bar z^2\\
\frac{1}{4}\tilde{g}=q_2dz^2+(|\mu|^2|\nu|^2\ti{h}^2+\ti{h}^{-2})dz\cdot d\bar z+\bar q_2d\bar z^2
\end{\eq}
So \begin{\eq}
\frac{1}{4}g(\frac{\partial}{\partial z},\frac{\partial}{\partial \bar{z}})=(|\nu|h-|\mu|h^{-1})^2+2|\mu\nu|,\quad \frac{1}{4}\tilde{g}(\frac{\partial}{\partial z},\frac{\partial}{\partial \bar{z}})=(|\mu||\nu|\ti{h}-\ti{h}^{-1})^2+2|\mu\nu|.
\end{\eq}
From Proposition \ref{harmonicmetriccomparison}, $|\mu||\nu|\ti{h}\leq |\mu|h^{-1}<\ti{h}^{-1}$, $|\mu||\nu|\ti{h}\leq |\nu|h<\ti{h}^{-1}$. Then
\begin{\eq}
(|\nu|h-|\mu|h^{-1})^2<(|\mu||\nu|\ti{h}-\ti{h}^{-1})^2,
\end{\eq}
which implies $g<\ti{g}$.\\

For $n=3$, we claim $|\nu|^2\ti{h}h^2<1$. The Hitchin equation is reduced to
\begin{\eq}
\triangle\log (|\nu|^2\ti{h}h^2)+\ti{h}^{-1}-|\mu|^4|\nu|^2\ti{h}^2+2|\mu|^2h^{-1}-2|\nu|^2h^2=0.
\end{\eq}
Let $u=|\nu|^2\ti{h}h^2$, $a=|\mu|^2\ti{h}h^{-1}$. Then
\begin{\eq}
\triangle\log u+\ti{h}^{-1}(1+2a-(2+a^2)u)=0.
\end{\eq}
Notice that $u\equiv1$ is a supersolution, then by the maximum principle, $u<1$. For the pullback metric $g,\ti{g}$, locally,
\begin{\eq}
\frac{1}{6}g=|\nu|^2h^2+2|\mu|^2h^{-1},\quad \frac{1}{6}\ti{g}=|\nu|^2|\mu|^4\ti{h}^2+2\ti{h}^{-1}.
\end{\eq}
Let $x=|\nu|h$, $\ti{x}=|\mu||\nu|^2\ti{h}$, $A=|q_3|=|\nu||\mu|^2$. Outside the zeros of $\mu\nu$, from Proposition \ref{harmonicmetriccomparison}, $x<\ti{x}$. Then
\begin{\eq}
\frac{1}{6}(g-\ti{g})&=&(x^2+\frac{2A}{x})-(\ti{x}^2+\frac{2A}{\ti{x}})
=\frac{(x-\ti{x})}{x\ti{x}}((x+\ti{x})x\ti{x}-2A)\\
&<&\frac{(x-\ti{x})}{x\ti{x}}(2x^2\ti{x}-2A)
=\frac{2|\nu||\mu|^2(x-\ti{x})}{x\ti{x}}(|\nu|^2\ti{h}h^2-1)<0.
\end{\eq}
So outside the zeros of $q_3=\mu^2\nu$, we obtain $g<\ti{g}$. We can easily see it also holds at the zeros of $q_3$.\\

For $n=4$, locally
\begin{\eq}
\frac{1}{8}g&=&|\nu|^2h^2_1+2|\ga|^2h_1^{-1}h_2+|\mu|^2h^{-2}_2
=(|\nu|h_1-|\mu|h^{-1}_2)^2+2|\nu||\mu|h_1h^{-1}_2+2|\ga|^2h_1^{-1}h_2,\\
\frac{1}{8}\ti{g}&=&|\mu|^2|\nu|^2|\ga|^4\ti{h}^2_1+2\ti{h}_1^{-1}\ti{h}_2+\ti{h}^{-2}_2
=(|\mu||\nu||\ga|^2\ti{h}_1-\ti{h}^{-1}_2)^2+2|\nu||\mu||\ga|^2\ti{h}_1\ti{h}^{-1}_2
+2\ti{h}_1^{-1}\ti{h}_2.
\end{\eq}
From Proposition \ref{harmonicmetriccomparison}, $|\mu||\nu||\ga|^{2}\ti{h}_1\leq |\mu|h^{-1}_2<\ti{h}^{-1}_2$, $|\mu||\nu||\ga|^{2}\ti{h}_1\leq |\nu|h_1<\ti{h}^{-1}_2$. Then
\begin{\eq}
(|\mu||\nu||\ga|^2\ti{h}_1-\ti{h}^{-1}_2)^2>(|\nu|h_1-|\mu|h^{-1}_2)^2.
\end{\eq}
Let $x=|\ga|^2h_1^{-1}h_2$, $\ti{x}=\ti{h}^{-1}_1\ti{h}_2$, $A=|q_4|=|\mu||\nu||\ga|^2$. \\

Claim:
$x<\ti{x}$ and $x\ti{x}> A$,
outside the zeros of $\ga$. Then the desired result follows from the basic identity $x+\frac{A}{x}-\ti{x}-\frac{A}{\ti{x}}=(x-\ti{x})(1-\frac{A}{x\ti{x}})$.

To show $x<\ti{x}$, let $u=\frac{x}{\ti{x}}=|\ga|^2h_1^{-1}h_2\ti{h}_1\ti{h}^{-1}_2$. Then $u$ satisfies
\begin{\eq}
\triangle \log u-2\ti{h}^{-1}_1\ti{h}_{2}(u-1)+2|\mu||\nu||\ga|^2\ti{h}_1\ti{h}^{-1}_{2}(u^{-1}-1)
+(|\nu|h_1-|\mu|h^{-1}_2)^2-(|\mu||\nu||\ga|^2\ti{h}_1-\ti{h}^{-1}_2)^2=0.
\end{\eq}
Then
\begin{\eq}
\triangle \log u-2\ti{h}^{-1}_1\ti{h}_{2}(u-1)+2|\nu||\mu||\ga|^2\ti{h}_1\ti{h}^{-1}_{2}(u^{-1}-1)>0.
\end{\eq}
Notice that $1$ is a subsolution, then by the maximum principle, $u<1$.

To show $x\ti{x}> A$, let $u=\frac{A}{x\ti{x}}=|\mu||\nu|h_1h^{-1}_2\ti{h}_1\ti{h}^{-1}_2$. Then $u$ satisfies
\begin{\eq}
\triangle \log u+(2\ti{h}^{-1}_1\ti{h}_2+2|\ga|^2h^{-1}_1h_2)(1-u)-(|\nu|h_1-|\mu|h^{-1}_2)^2-(|\mu||\nu||\ga|^2\ti{h}_1-\ti{h}^{-1}_2)^2=0.
\end{\eq}
Then
\begin{\eq}
\triangle \log u+(2\ti{h}^{-1}_1\ti{h}_2+2|\ga|^2h^{-1}_1h_2)(1-u)>0.
\end{\eq}
Notice that $u\equiv1$ is a solution, then by the maximum principle, $u<1$. At the zeros of $\ga$, we can also obtain $g<\ti{g}$ from $|\mu||\nu|h_1h^{-1}_2\ti{h}_1\ti{h}^{-1}_2<1$. So we finish the proof.
\qed

By integration, we obtain
\begin{cor}
The Morse function achieves the maximum in the Hitchin point in the above cases.
\end{cor}
As an immediate corollary in terms of representations for $n=2$, we recover the following result shown in \cite{DominationFuchsian}.
\begin{cor} For any non-Fuchsian reductive $SL(2,\mathbb{R})$-representation $\rho$ and any Riemann surface $\Sigma$, there exists a Fuchsian representation $j$ such that the pullback metric of the corresponding $j$-equivariant harmonic map $f_j:\widetilde{\Sigma}\rightarrow \mathbb{H}^2$ dominates the one for $f_{\rho}$. \end{cor}
\begin{pf}
For any reductive $SL(2,\mathbb{R})$-representation $\rho$, if it is into the compact subgroup $SO(2,\mathbb{R})$, the associated harmonic map is constant. In this case, the statement is clear. Given any Riemann surface $\Sigma$, if the representation $\rho$ is not into the compact group $SO(2,\mathbb{R})$, it corresponds to a cyclic Higgs bundle parametrized by $(\alpha,\beta)$ over $\Sigma$ by \cite{Hitchin87}. Then we choose the Fuchsian representation $j$ corresponding to the cyclic Higgs bundle parametrized by $q_2=\alpha\beta$ over $\Sigma$. The statement follows from Theorem \ref{comparisonarticle}.
\end{pf}\qed

\section{Maximal $Sp(4,\mathbb{R})$-representations}\label{maximal}
For each reductive representation $\rho$ into $Sp(2n,\mathbb{R})$, we can define a Toledo integer $\tau(\rho):=\frac{2}{\pi}\int_Sf^*\omega$ where $f$ is any $\rho$-equivariant continuous map $f:\widetilde{S}\rightarrow Sp(2n,\mathbb{R})/U(n)$ and $\omega$ is the normalized $Sp(2n,\mathbb{R})$-invariant K\"ahler $2$-form on $Sp(2n,\mathbb{R})/U(n)$. It is well-known that $|\tau(\rho)|\leq n(g-1)$. The representation $\rho$ with $|\tau(\rho)|=n(g-1)$ is called maximal.

A $Sp(4,\mathbb{R})$-Higgs bundle over $\Sigma$ is a pair $(V\oplus V^*,\begin{pmatrix}0&\beta\\ \gamma&0\end{pmatrix})$ where $V$ is a rank $2$ holomorphic vector bundle over $\Sigma$, $\beta\in H^0(S^2V\otimes K_{\Sigma})$ and $\gamma\in H^0(\Sigma,S^2V^*\otimes K_{\Sigma})$. The Toledo integer of the $Sp(4, \mathbb{R})$-Higgs bundle is the integer $\deg(V)$. 
There are $3 \cdot 2^{2g} + 2g -4$ components of maximal $Sp(4,\mathbb{R})$-representations shown in \cite{Gothen} containing $2^{2g}$ Hitchin components isomorphic to each other and $2g-3$ exceptional components called Gothen components.

Labourie in \cite{LabourieCyclic} shows that any $Sp(4,\mathbb{R})$ Hitchin representation corresponds to a cyclic Higgs bundle in the Hitchin components over a unique Riemann surface. As a result, there is a unique $\rho$-equivariant minimal immersion of $\widetilde{S}$ into $Sp(4,\mathbb{R})/U(2)$ for any Hitchin representation for $Sp(4,\mathbb{R})$.

For each Riemann surface $\Sigma$, each Gothen component is explicitly described in \cite{BradlowDeformation} as the moduli space of Higgs bundles of the following form
\begin{\eq}
E=N\oplus NK^{-1}\oplus N^{-1}K\oplus N^{-1}, \quad \phi=
\left(
\begin{array}{cccccc}
0 &q_2 &0 & \nu\\
1 &0 &0 &0 \\
0&\mu &0 &q_2 \\
0 & 0& 1 & 0
\end{array}
\right)
\end{\eq} where $g-1<\deg(N)<3g-3$, $\nu\in H^0(N^2K), \mu\in H^0(N^{-2}K^3),$ and $q_2\in H^0(K^2)$. Here $V=N\oplus N^{-1}K$.

By Collier's work \cite{Collier}, we can replace the variation of $q_2$ with a variation of base Riemann surface structure. That is, any maximal $Sp(4,\mathbb{R})$-representation in the Gothen components corresponds to a Higgs bundle over a unique Riemann surface $\Sigma$ of the form
\begin{\eq}
E=N\oplus NK^{-1}\oplus N^{-1}K\oplus N^{-1}, \quad \phi=
\left(
\begin{array}{cccccc}
0 &0 &0 & \nu\\
1 &0 &0 &0 \\
0&\mu &0 &0 \\
0 & 0& 1 & 0
\end{array}
\right)
\end{\eq} where $g-1<\deg(N)<3g-3$, $\mu\neq 0$ and $\nu$ can be zero. The $2g-3$ Gothen components are indexed by the degree of $N$. These are cyclic $SL(4,\mathbb{C})$-Higgs bundles. Note that If $N=K^{\frac{3}{2}}$, this gives the Hitchin representation. As a result, for any $Sp(4,\mathbb{R})$-representation in the Gothen components, there is a unique $\rho$-equivariant minimal immersion of $\widetilde{S}$ into $Sp(4,\mathbb{R})/U(2)$.

The above cyclic Higgs bundles with $\nu=0$ are stable and play a similar role as the Fuchsian case. We call the corresponding representations $\mu$-Fuchsian representations. The space of $\mu$-Fuchsian representations serves as the minimum in its component of maximal $Sp(4,\mathbb{R})$ representations in the following sense.
\begin{cor}
For any maximal representation $\rho:\pi_1(S)\rightarrow Sp(4,\mathbb{R})$ in the $2g-3$ Gothen components, there exists a $\mu$-Fuchsian representation $j$ of $\pi_1(S)$ such that the pullback metric of the unique $j$-equivariant minimal immersion $f_j:\widetilde{S}\rightarrow Sp(4,\mathbb{R})/U(2)$ is dominated by the one for $f_{\rho}$.
\end{cor}
\begin{pf}
For any maximal representation in the Gothen component, we can realize it as a cyclic Higgs bundle parametrized by $(1,\mu,1,\nu)$ over some Riemann surface $\Sigma$. Then we choose the $\mu$-Fuchian representation corresponding to cyclic Higgs bundle parametrized by $(1,\mu,1,0)$ over $\Sigma$. Then the statement follows from Theorem \ref{domination2article}.
\end{pf}\qed\\

Since any Hitchin representation for $Sp(4,\mathbb{R})$ corresponds to a cyclic Higgs bundle over some Riemann surface $\Sigma$, we obtain bounds on the extrinsic curvature of minimal immersions for maximal representations in the Hitchin component as an immediate corollary of Theorem \ref{curvaturearticle}.
\begin{cor} For any Hitchin representation $\rho$ for $Sp(4,\mathbb{R})$, the sectional curvature $K_{\sigma}$ in $Sp(4,\mathbb{R})/U(2)$ of the tangent plane $\sigma$ of the uniuqe $\rho$-equivariant minimal immersion satisfies \\
(1) $K_{\sigma}=-\frac{1}{40}$, if $\rho$ is Fuchsian;\\
(2) $-\frac{1}{36}<K_{\sigma}<0$ and $\exists~ p$ such that $K_{\sigma}(p)<-\frac{1}{40},$ if $\rho$ is not Fuchsian.
\end{cor}
\begin{rem} The lower bound $-\frac{1}{36}$ is nontrivial, since the sectional curvature $K$ in $Sp(4,\mathbb{R})/U(2)$ satifies that $-\frac{1}{4}\leq K\leq 0$.
\end{rem}
Similarly, we also obtain estimates on the extrinsic curvature of minimal immersions for maximal representations in $2g-3$ Gothen component.
\begin{thm}
For any maximal representation $\rho$ for $Sp(4,\mathbb{R})$ in each Gothen component, the sectional curvature $K_{\sigma}$ in $Sp(4,\mathbb{R})/U(2)$ of the tangent plane $\sigma$ of the unique $\rho$-equivariant minimal immersion satisfies \\
(1) $-\frac{1}{8}\leq K_{\sigma}<-\frac{1}{40}$ and the lower bound is sharp, if $\rho$ is $\mu$-Fuchsian;\\
(2) $-\frac{1}{8}\leq K_{\sigma}<0$, if $\rho$ is not $\mu$-Fuchsian.
\end{thm}

\begin{pf} It is sufficient to work with cyclic Higgs bundle parameterized by $(1,\mu,1,\nu)$ of the above form.
The Hitchin equation in this case is
\begin{\eq}
\triangle \log {h}_1+{h}_1^{-1}h_{2}-|\nu|^2{h}_{1}^{2}&=&0,\\
\triangle \log {h}_2+|\mu|^2{h}_2^{-2}-h_1^{-1}h_2&=&0
\end{\eq}

Using the curvature formula (\ref{CurvatureFormula}), the sectional curvature of the tangent plane $\sigma$ of the minimal immersion is
\begin{\eq}
K_{\sigma}&=&-\frac{(h_1^{2}|\nu|^2-h_1^{-1}h_{2})^2+(h_1^{-1}h_2-h_2^{-2}|\,u|^2)^2}{4\cdot(h_1^{2}|\nu|^2+2h_1^{-1}h_{2}+h_{2}^{-2}|\mu|^2)^2}.
\end{\eq}
For the right inequality, outside zeros of $\mu\nu$,
\begin{\eq}
&&\triangle \log {h}_1^2h_2^{-2}|\mu\nu|-(|\mu|^2h_2^{-2}+|\nu|^2{h}_{1}^{2})+2{h}_1^{-1}h_{2}=0\\
\Longrightarrow&& \triangle \log {h}_1^2h_2^{-2}|\mu\nu|-2|\mu\nu|h_1h_2^{-1}-2{h}_1^{-1}h_{2}\geq 0\\
\Longrightarrow&& \triangle \log {h}_1^2h_2^{-2}|\mu\nu|-2(h_1^2h_2^{-2}|\mu\nu|-1){h}_1^{-1}h_{2}\geq 0
\end{\eq}
So at the maximum of $h_1^2h_2^{-2}|\mu\nu|$, $h_1^2h_2^{-2}|\mu\nu|-1\leq 0$. Hence $h_1^2h_2^{-2}|\mu\nu|\leq 1$ on the whole surface. By the strong maximum principle, we obtain that $h_1^2h_2^{-2}|\mu\nu|<1$. So ${h}_1^{-1}h_{2}=|\nu|^2{h}_{1}^{2}$ and $|\mu|^2{h}_2^{-2}=h_1^{-1}h_2$ cannot hold at any point $p$ simultaneously, since it would imply that $h_1^2h_2^{-2}|\mu\nu|=1$ at point $p$, contradiction. Therefore $K_{\sigma}<0$.\\\\
For the left inequality. Let $f_1=\frac{h_1^2|\nu|^2}{h_1^{-1}h_2}, f_2=\frac{h_2^{-2}|\mu|^2}{h_1^{-1}h_2}$. Claim: $f_1,f_2<\frac{4}{3}$.
\\
The equation for $f_1$ is, outside zeros of $\nu$,
\begin{\eq}
&&\triangle \log f_1+[3(1-f_1)-(f_2-1)]h_1^{-1}h_2=0\\
\Longrightarrow &&\triangle_{h_1^{-1}h_2}\log f_1+3(1-f_1)-(f_2-1)=0\\
\Longrightarrow &&\triangle_{h_1^{-1}h_2}\log f_1+3(1-f_1)+1\geq 0
\end{\eq}
So at the maximum of $f_1$, $3(1-f_1)+1\leq 0$, hence $f_1\leq \frac{4}{3}$. Use the strong maximum principle, $f_1<\frac{4}{3}$.
It is similar for $f_2$. The claim is proven.

Using $0\leq f_1,f_2<\frac{4}{3}$,
\begin{\eq}
K_{\sigma}=-\frac{\text{tr}([\phi,\phi^*][\phi,\phi^*])}{2n\cdot\text{tr}(\phi\phi^*)^2}=-\frac{(f_1-1)^2+(f_2-1)^2}{4(2+f_1+f_2)^2}\geq-\frac{1+1}{16}=-\frac{1}{8}.
\end{\eq}

Note that $K_{\sigma}$ only achieves $-\frac{1}{8}$ if $f_1=f_2=0$. This only happens at common zeros of $\mu$ and $\nu$.

In the $\mu$-Fuchsian case, $\nu=0$. So $f_1=0$ and again $f_2<\frac{4}{3}$. Then using $(f_2+2)^{-1}\in (\frac{3}{10},\frac{1}{2}]$,
\begin{\eq}
K_{\sigma}= -\frac{1+(f_2-1)^2}{4(f_2+2)^2}=-\frac{(f_2+2-3)^2+1}{4(f_2+2)^2}=-\frac{10}{4}(((f_2+2)^{-1}-\frac{3}{10})^2-\frac{1}{40}<-\frac{1}{40}.
\end{\eq}
Note that at zeros of $\mu$ in $\mu$-Fuchsian case, the curvature $K_{\sigma}=-\frac{1}{8}$.
\end{pf}\qed
\begin{rem}
As shown in \cite{Mochizuki}, along the family of $(E,t\phi)$, as $|t|\rightarrow \infty$, away from zeros of $\mu\nu\neq 0$, the sectional curvature goes to zero. So the upper bound in Part (2) is sharp.
\end{rem}

We compare the Gothen components with the Hitchin components.
\begin{cor}
For any maximal representation $\rho:\pi_1(S)\rightarrow Sp(4,\mathbb{R})$ in the $2g-3$ Gothen components, there exists a Hitchin representation $j$ of $\pi_1(S)$ such that the pullback metric of the unique $j$-equivariant minimal immersion $f_j:\widetilde{S}\rightarrow Sp(4,\mathbb{R})/U(2)$ dominates the one for $f_{\rho}$.
\end{cor}
\begin{pf}
For any maximal representation $\rho$ in the Gothen components, it corresponds to a cyclic Higgs bundle parametrized by $(1,\mu,1,\nu)$ over some Riemann surface $\Sigma$. Then we choose the Hitchin representation $j$ corresponding to cyclic Higgs bundle in the Hitchin component parametrized by $q_4=\mu\nu$ over $\Sigma$. The statement then follows from Theorem \ref{comparisonarticle} for $n=4$.
\end{pf}\qed

\end{document}